\documentclass[preprint,12pt]{elsarticle}

\usepackage{amssymb,amsmath}
\usepackage{epsfig, graphicx}
\graphicspath{ {./figures/} }
\usepackage[all]{xy}
\usepackage{color}
\usepackage{enumitem}  
\usepackage{ulem} 
\usepackage{mathrsfs} 
\usepackage{multirow}
\usepackage{subcaption} 
\usepackage{tikz}
\usetikzlibrary{shapes, arrows, fit, positioning, arrows.meta}


\newcommand{\RR}{\mathbb{R}}
\newcommand{\xx}{{\mathbf{x}}}
\newcommand{\ttheta}{\boldsymbol{\theta}}

\newtheorem{remark}{Remark}

\newtheorem{assumption}{Theorem}
\newtheorem{assumptionalt}{Theorem}[assumption]
\newenvironment{assumptionp}[1]
{
	\renewcommand\theassumptionalt{#1}
	\assumptionalt
}{\endassumptionalt}

\journal{Journal of Computational and Applied Mathematics \hspace{1em}}

\begin{document}
\begin{frontmatter}
\title{A Machine Learning and Finite Element Framework for Inverse Elliptic PDEs via Dirichlet-to-Neumann Mapping}
\author[KNU1]{Dabin Park}
\author[FSU]{Sanghyun Lee}
\author[KNU2]{Sunghwan Moon}
\affiliation[KNU1]{organization={School  of  Mathematics,  Kyungpook  National  University},  city={Daegu},  postcode={41566},  country={Republic  of  Korea}}
\affiliation[FSU]{organization={Department of Mathematics, Florida State University}, city={Tallahassee}, state={FL}, postcode={32306}, country={United States}}
\affiliation[KNU2]{organization={Department  of  Mathematics,  Kyungpook  National  University},  city={Daegu},  postcode={41566},  country={Republic  of  Korea}}

\begin{abstract}
Inverse problems for Partial Differential Equations (PDEs) are crucial in numerous applications such as geophysics, biomedical imaging,
and material science, where unknown physical properties must be inferred from indirect measurements.
In this work, we present a new approach to solving the inverse problem for elliptic PDEs, using only boundary data.
Our method leverages the Dirichlet-to-Neumann (DtN) map, which captures the relationship between boundary inputs and flux responses.
This enables the reconstruction of the unknown physical properties within the domain from boundary measurements alone.
Our framework employs a self-supervised machine learning algorithm that integrates a Finite Element Method (FEM) in the inner loop for the forward problem, ensuring high accuracy.
Moreover, our approach illustrates its effectiveness in challenging scenarios with only partial boundary observations, which is often the case in real-world scenarios.
In addition, the proposed algorithm effectively handles discontinuities by incorporating carefully designed loss functions.
This combined FEM and machine learning approach offers a robust, accurate solution strategy for a broad range of inverse problems,
enabling improved estimation of critical parameters in applications from medical diagnostics to subsurface exploration.
\end{abstract}


\end{frontmatter}

\section{Introduction}
\label{sec:intro}

 Inverse problems involve inferring unknown properties from indirect observations, in contrast to forward problems, where the governing equations and parameters are given and the solution is computed directly.
Such problems often arise in science and engineering, particularly when reconstructing coefficients, sources, or initial conditions in governing equations.
In particular, inverse problems for Partial Differential Equations (PDEs) are a fundamental area of study in applied and computational mathematics~\cite{Isakov06, Kirsch11},
with wide-ranging applications in fields such as geophysics~\cite{Duchateau97}, biomedical imaging~\cite{Kuchment13, Steinberg19}, and material science~\cite{Janno11,Lee25}.
For instance, in medical imaging techniques such as Electrical Impedance Tomography (EIT), recovering these hidden material properties from indirect measurements enables visualization of conductivity distributions inside the human body, thereby enhancing diagnostic capabilities~\cite{Bar21, Cheney99, Patterson05}.
Similarly, in geophysical contexts, inferring subsurface properties such as permeability from surface measurements is critical for resource exploration and environmental monitoring~\cite{Ewing91, Carrera88}.

In recent years, machine learning techniques have been incorporated into the study of inverse problems, offering new perspectives beyond traditional approaches.
Pioneering works have demonstrated that deep neural networks can solve PDE-related inverse problems directly, and improve reconstruction compared to classical methods~\cite{Jin17, Adler18}.
These developments have been systematically reviewed and categorized in comprehensive surveys on deep learning for inverse problems~\cite{Arridge19, Ongie20}.
Building on these advances, the neural operator framework was recently introduced to approximate solution operators between function spaces~\cite{Kovachki23}, enabling the efficient resolution of parameterized PDE inverse problems.
In parallel, the Finite Element Methods (FEM) has long played a central role in inverse problems, particularly in elasticity imaging.
Early studies used the FEM to reconstruct material parameters from displacement data, and adjoint-based formulations later improved efficiency for biomedical applications~\cite{Maniatty89,Oberai03,Doyley12}.
FEM has also been applied to other inverse settings such as optical tomography~\cite{Arridge09}, and recent work has extended its use to nonlinear constitutive models and hybrid approaches that integrate FEM with neural networks for efficient PDE inversion~\cite{Xu25,Meethal23}.

\begin{figure}[!ht]
    \centering
    \begin{tikzpicture}[scale=0.75, font=\small]
    
    \draw[thick,fill=yellow!10] (-2,-2) rectangle (2,2);
    \node[] at (0,2.4) {Conductivity ($k_i$) is known};
    
    \node[] at (8,2.4) {Obtain the streamline};

    \draw[->,thick] (2.5,1.5) -- (5.5,1.5) 
        node[midway,above] {Forward problem};
    \draw[->,thick] (5.5,-1.5) -- (2.5,-1.5)
        node[midway,below] {Inverse problem};

    \draw[fill=yellow!50,draw=black] (-0.8,0.6) circle (0.4)     node {\(k_1\)};
    
    \draw[rotate=45,fill=blue!20,draw=black] (0.7,-0.5) ellipse (10pt and 20pt)
    node {\(k_2\)};
          
    \node[] at (-0.2,-0.9) {\(k_3\)};
    
    \draw[thick, blue, dash pattern=on 8pt off 8pt] (6,-2) rectangle (10,2);
    
    \node[] at (8,-2.4) {Boundary measurement is known};
    
    \draw[thick] (10, 1.5) .. controls (8, 1.7) .. (6, 1.5);
    \draw[thick] (10, 1.0) .. controls (8, 1.2) and (7, 1.3) .. (6, 1.0);
    \draw[thick] (10, 0.7) .. controls (9, 0.5) and (7, 0.9) .. (6, 0.7);
    \draw[thick] (10, 0.0) .. controls (8, 0.0) .. (6, 0.0);
    \draw[thick] (10, -0.7) .. controls (9, -0.9) and (7, -0.3) .. (6, -0.7);
    \draw[thick] (10, -1.0) .. controls (8, -1.2) and (7, -1.3) .. (6, -1.0);
    \draw[thick] (10, -1.5) .. controls (8, -1.7) .. (6, -1.5);
    
    \end{tikzpicture}
    \caption{Illustration of the problem. The dashed line indicates the partially observed boundary measurement.}
    \label{fig:intro}
\end{figure}
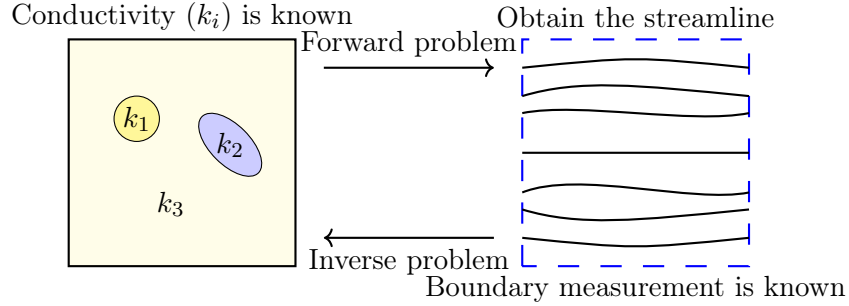

In this paper, we focus on the inverse problem for elliptic PDEs with variable, potentially discontinuous coefficients.
The main novelty of our approach is the utilization of the Dirichlet-to-Neumann (DtN) map—a mathematical construct that encapsulates the relationship between prescribed boundary inputs and their resulting flux responses~\cite{Feldman19}.
This mapping framework is especially advantageous because it allows for the reconstruction of the interior characteristics of a domain using solely boundary measurements.
In scenarios such as EIT or subsurface exploration, where sensors are predominantly deployed on the boundary, this capability eliminates the need for invasive internal data collection while still providing detailed insight into the material properties within the domain.
See Figure~\ref{fig:intro} for the illustration.

In addition, our approach is flexible to utilize partial observations of the boundary data, which is often the case in real-world scenarios. The presented numerical examples further illustrate the capabilities and robustness of the proposed algorithm, demonstrating its effectiveness across a range of challenging inverse problems.

The proposed computational framework further distinguishes itself by integrating a self-supervised machine learning algorithm~\cite{Subramanian23} with FEM used in the inner loop for solving the forward problem.
This combination not only achieves high accuracy in modeling the physical processes but also enables iterative refinement of the coefficient reconstructions. The carefully designed loss functions, tailored to handle discontinuities effectively, enhance the overall fidelity of the inverse solution, making our approach robust and versatile. Ultimately, this framework offers a precise and reliable estimation of critical parameters, paving the way for significant advancements in applications ranging from non-invasive medical diagnostics to environmental and subsurface monitoring.

In summary, our novelties include \begin{enumerate}
    \item Reconstruction of an unknown coefficient \textbf{using only boundary measurements}, without relying on interior data.
    \item The framework demonstrates good performance even with \textbf{partial boundary observations}, reflecting realistic sensing scenarios.
    \item It can reconstruct \textbf{discontinuous coefficients} by incorporating carefully designed loss functions.
    \item A \textbf{self-supervised framework} is employed that learns from Cauchy data pairs, without requiring labeled interior coefficients.
     \item The proposed method \textbf{integrates the Finite Element Method (FEM) into the inner loop with adjoint-based gradients}, utilizing FEM as a physics-grounded solver to ensure accuracy and consistency with the underlying PDE while maintaining computational efficiency in the inversion process.
\end{enumerate}

The rest of the paper is organized as follows. First, the governing system and the DtN mapping are introduced in Section \ref{sec:2}.
The overall computational framework, including the loss functions, is discussed in Section \ref{sec:algorithm}.
Finally, in Section \ref{sec:num}, several numerical experiments to illustrate the capabilities and the performance of our algorithm are presented.

\section{Governing System}
\label{sec:2}
In this paper, we consider the following elliptic PDE
\begin{equation}\label{eq:Elliptic}
\left\{\begin{array}{rll}
-\nabla \cdot (k(\xx) \nabla p(\xx)) &= f(\xx)  &\xx \in \Omega ,\\
p(\xx) &= g(\xx)  &\xx \in \partial \Omega ,
\end{array}\right.
\end{equation}
where $\Omega$ is a bounded domain in $\RR^2$,
$p: \overline{\Omega} \rightarrow \RR$ is the solution function,
$f$ is the source function, $g \in H^{1/2}(\partial \Omega)$ is the given Dirichlet boundary condition, 
and $k$ is a strictly positive coefficient defined in $\overline{\Omega}$.
Here, $\xx = (x,y)$ denotes the spatial coordinate in $\RR^2$.

We consider an inverse problem associated with \eqref{eq:Elliptic}, where the goal is to reconstruct the coefficient $k$ from the boundary measurements.
To achieve this, we employ deep learning techniques.
A straightforward approach would be to use supervised learning, which requires access to exact coefficient values of $k(\xx)$ for all $\xx \in \overline{\Omega}$.
However, in most cases, obtaining complete knowledge of $k$ in $\overline{\Omega}$ is infeasible.
Thus, we propose a self-supervised learning framework that does not rely on exact coefficient values of $k$ but instead leverages available data such as Dirichlet and Neumann boundary conditions to infer the coefficient.

The methodology and key concepts necessary for solving this problem are introduced in the following subsection.

\subsection{Dirichlet-to-Neumann Map}
\label{sec:dtn}
In this subsection, we introduce the DtN map, which is central to the reconstruction of $k$.
This map is defined as follows:
\begin{equation*}
\Lambda_k : H^{1/2}(\partial\Omega) \to H^{-1/2}(\partial\Omega),
\end{equation*}
which maps the Dirichlet data $g$ to the corresponding Neumann data $h_g = k \partial_\nu p_g|_{\partial \Omega}$.
Here, $p_g$ is the solution of \eqref{eq:Elliptic} for the given source function $f$ and Dirichlet condition $g$,
and $\nu$ denotes the unit outward normal vector on $\partial \Omega$.
The space $H^r(\partial\Omega)$, $r \in \RR$, denotes the Sobolev space.

The reconstruction of the coefficient $k$ using the DtN map has been extensively studied in the context of EIT and the Calder\'on problem~\cite{Feldman19}.
In particular, the case when $f=0$ has been comprehensively surveyed, including the uniqueness of $k$ under partially observed boundary data~\cite{Uhlmann14}.

Here, we briefly recapitulate the theorems that guarantee the existence and the uniqueness of $k$ with given conditions.
Let $\Gamma \subset \partial \Omega$ be a non-empty open subset of the boundary. Then the following theorem proves the uniqueness of $k$ in the case where $f =0$.
\begin{assumptionp}{A}[Corollary 1.19 in~\cite{Uhlmann14}]
    Let $\Omega \subset \RR^2$ be a bounded domain with a smooth boundary, and let $k_1, k_2 \in C^{4+\alpha}(\overline{\Omega})$ be non-vanishing functions for some $\alpha > 0$.
    Assume that  
    \begin{equation*}
    \Lambda_{k_1}(g) = \Lambda_{k_2}(g) \text{ on } \Gamma \text{ for all } g \in H^{1/2}(\partial\Omega),\, \operatorname{ supp} g \subset \Gamma.
    \end{equation*}  
    Then, $k_1 = k_2$ in $\Omega$.
\end{assumptionp}
Furthermore, if $\Gamma = \partial \Omega$, the uniqueness of $k$ can be established using data from the entire boundary.

Next, in the case where $f \neq 0$, we note that the uniqueness of $k$ in the entire domain $\Omega$ is not guaranteed. However, its uniqueness on the boundary $\partial \Omega$ has been established in~\cite{Isakov06} as the following theorem:
\begin{assumptionp}{B}[Theorem 5.1.1 in~\cite{Isakov06}]
    Assume that $k_1,k_2 \in C^m(V)$, where $V$ is a neighborhood of a boundary point of $\Omega$.
    \begin{equation*}
        \text{ If } \Lambda_{k_1} = \Lambda_{k_2}, \text{ then } \partial^\alpha k_1 = \partial^\alpha k_2 \text{ on } \partial\Omega\cap V \text{ when } |\alpha| \leq m.
    \end{equation*}
\end{assumptionp}
If $k_1,k_2\in C(\Omega)$ are two coefficients that are analytic in $\Omega$ and $\Lambda_{k_1} = \Lambda_{k_2}$, then $k_1=k_2$ on $\partial \Omega$ and by Strong Maximum Principle \cite[Chapter 2]{evans10}, we have $k_1 = k_2$ in $\Omega$.

Let $\omega$ be either $\Gamma$ or $\partial \Omega$.
We refer to the case $\omega = \Gamma$ as the partially observed setting, and $\omega = \partial \Omega$ as the fully observed setting.
Instead of directly using the DtN map, which is usually difficult to obtain in real-world scenarios, we utilize Cauchy data—pairs of Dirichlet and Neumann boundary measurements.
\begin{equation*}
\left\{(g, h_g) : g \in H^{1/2}(\partial\Omega), \, \operatorname{supp} g \subset \omega, h_g = k \partial_\nu p_g|_{\omega} \right\}.
\end{equation*}
These pairs serve as practical substitutes for the DtN map in our setup.
In practice, it is not feasible to train on all functions $g \in H^{1/2}(\partial\Omega)$ with supp~$ g \subset \omega$.
We consider the following problem:\\[0.2cm]
\noindent\textbf{Problem.} Given a sufficiently large finite collection of Cauchy data:
\begin{equation*}
    \left\{ (g_i, h_{g_i})_{i = 1, \dots, N}: g_{i} \in H^{1/2}(\partial\Omega), \, \operatorname{supp} g_{i} \subset \omega, h_{g_i} = k \partial_\nu p_{g_i}|_{\omega} \right\},
\end{equation*}
the coefficient $k$ can be numerically reconstructed in $\Omega$.

\begin{remark}{}
The above theorems establish the theoretical uniqueness of the coefficient $k$ under certain conditions on the DtN map.
Based on these results, our work reconstructs $k$ from finitely many Cauchy data, as utilizing the full DtN map is not feasible in practice.
Moreover, our framework handles both the fully observed setting ($\omega = \partial\Omega$) and the partially observed setting ($\omega = \Gamma$).
\end{remark}

Rather than directly using the values of $k$, we adopt a self-supervised learning framework where the model learns to reconstruct $k$ from Cauchy data alone~\cite{Hwang23, Hwang24}.
This approach is particularly useful in practical setups where $k$ is not directly measurable.
\section{Computational Algorithm}
\label{sec:algorithm}

In this section, we present the numerical approach for solving \textbf{Problem}.
First, in Subsection~\ref{subsec:framework}, we provide an overview of the proposed framework.
Then, we introduce the model for reconstructing the coefficient $k$ in Subsection~\ref{subsec:model_k}.
Next, in Subsection~\ref{subsec:forward}, we briefly describe how to solve the forward problem using the FEM.
Finally, in  Subsection~\ref{subsec:loss}, we define the loss function used to train the model and describe the method for calculating its gradient to update the model parameters.

\subsection{Overall Framework}
\label{subsec:framework}

In this subsection, we introduce the proposed computational framework for reconstructing the coefficient $k$.

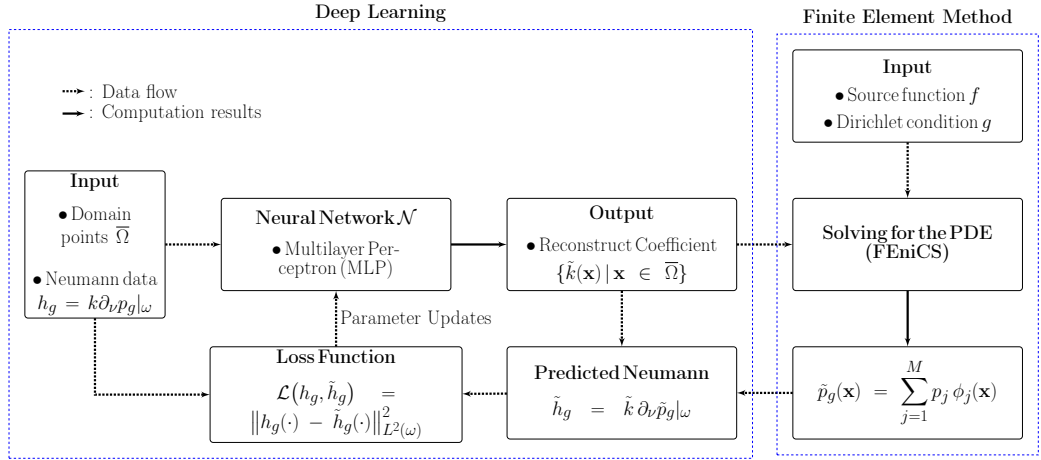
\begin{figure}[!ht]
    \centering
    \resizebox{\textwidth}{!}{
\begin{tikzpicture}[%
    auto,
    node distance=2.5cm,
    >=latex',
    every node/.style={font=\small}
]

\tikzset{
  block/.style={
    rectangle,
    draw,
    rounded corners,
    align=center,
    text width=5cm,
    minimum height=10em
  },
  bigblock/.style 2 args={
    draw=blue,
    line width=10mm, 
    thick,
    dashed,
    inner sep=1.7em,
    label={[anchor=north, yshift=3em]north:\textbf{\large #1}}, 
    fit=#2
  },
  line/.style={
    draw,->,
    >=latex'
  }
}

\node[block, text width=6cm] (input) {
  \textbf{\Huge Input}\\[1.0em]
  \huge $\bullet$ Domain points \;$\overline{\Omega}$\\[1.0em]
  \huge $\bullet$ Neumann data\\[0.5em]
  \huge $h_g = k \partial_\nu p_g|_\omega$
};

\node[block, right=of input, text width=10cm] (nn) {
  \textbf{\Huge Neural Network $\mathcal{N}$}\\[2.0em]
  \huge $\bullet$ Multilayer Perceptron (MLP)
};

\node[block, text width=3cm, above=0.5cm of input, draw=white, align=left] (legend) {
    \begin{tikzpicture}[baseline=(current bounding box.north), every node/.style={font=\Huge}]
        \draw[->, line width=1mm] (4,-2) -- (5,-2) node[right] {: Computation results};
        \draw[->, dotted, line width=1mm] (4,-1) -- (5,-1) node[right] {: Data flow};
    \end{tikzpicture}
};

\node[block, above=of nn, draw opacity=0.0] (transparent) {
};

\node[block, below=of nn, text width=11cm] (loss) {
  \textbf{\Huge Loss Function}\\[2.0em]
  \huge $\displaystyle
    \mathcal{L}\bigl(h_{g}, \tilde{h}_g \bigr)
    = \bigl\|h_{g}(\cdot) -\tilde{h}_g(\cdot)
    \bigr\|^2_{L^2(\omega)}
  $
};

\node[block, right=of nn, text width=10cm] (output) {
  \textbf{\Huge Output}\\[1.0em]
  \huge $\bullet$ Reconstruct Coefficient\\[0.5em]
  \Huge $\{\tilde{k}(\xx) \,|\, \xx \!\! \in \!\! \overline{\Omega} \}$
};

\node[block, below=of output, text width=10cm] (predicted_neumann) {
  \textbf{\Huge Predicted Neumann}\\[2.0em]
  \Huge $\tilde{h}_g = \tilde{k}\,\partial_\nu \tilde{p}_{g}|_\omega$
};

\node[block, right=of output, text width=10cm] (fenics) {
  \textbf{\Huge Solving for the PDE \\ (FEniCS)}
};

\node[block, above=of fenics, text width=10cm] (input2) {
  \textbf{\Huge Input}\\[0.5em]
  \huge $\bullet$ Source function $f$ \\[0.5em]
  \huge $\bullet$ Dirichlet condition $g$
};

\node[block, below=of fenics, text width=10cm] (predicted_sol) {
  \Huge $\displaystyle \tilde{p}_{g}(\xx)=\sum_{j=1}^M p_j\,\phi_j(\xx)$
};

\path[line, line width=1mm, dotted] (input) -- (nn);
\path[line, line width=1mm, dotted] (input.south) |- (loss.west);
\path[line, line width=1mm]         (nn) -- (output);
\path[line, line width=1mm, dotted] (output) -- (fenics);
\path[line, line width=1mm, dotted] (input2) -- (fenics);
\path[line, line width=1mm]         (fenics) -- (predicted_sol);
\path[line, line width=1mm, dotted] (output) -- (predicted_neumann);
\path[line, line width=1mm, dotted] (predicted_sol) -- (predicted_neumann);
\path[line, line width=1mm, dotted] (predicted_neumann) -- (loss);
\path[line, line width=1mm, dotted] (loss) -- (nn) node[midway, right] {\Huge Parameter Updates};

\node[bigblock={\Huge Deep Learning}{(legend) (input) (transparent) (nn) (loss) (output) (predicted_neumann) }] (dlbox) {};
\node[bigblock={\Huge Finite Element Method}{(input2) (fenics) (predicted_sol)}] (fembox) {};

\end{tikzpicture}
}

\caption{Framework for reconstructing $k$ using a neural network and FEM.}
\label{fig:framework}

\end{figure}

Figure~\ref{fig:framework} illustrates the process of reconstructing the coefficient $k$ using a neural network and FEM.
The framework proceeds as follows:
First, the domain points $\overline{\Omega}$ are the input to the neural network $\mathcal{N}$, which outputs the reconstruction $\tilde{k}$ of $k$.
The neural network $\mathcal{N}$ is designed as a Multilayer Perceptron (MLP) that reconstructs $k$ from the input domain points.
Next, the reconstructed coefficient $\tilde{k}$, along with the Dirichlet condition $g$ and source function $f$, is used to solve \eqref{eq:Elliptic}.
The resulting predicted solution is denoted as $\tilde{p}_g$, which is represented as a linear combination of basis functions $\{ \phi_i \}_{i=1}^M$ in the finite-dimensional subspace $V_h$.
These basis functions and the finite-dimensional subspace are further detailed in Subsection~\ref{subsec:forward}.
Using $\tilde{k}$ and $\tilde{p}_g$, the predicted Neumann data $\tilde{h}_g =\tilde{k}\partial_\nu \tilde{p}_g|_\omega$ is computed.
This predicted Neumann data is then compared with the given Neumann data $h_g$ to calculate the loss function.
Finally, the model parameters of the neural network $\mathcal{N}$ are updated to minimize the loss function.
By iterating through this process, the neural network $\mathcal{N}$ learns to reconstruct $k$ more accurately.

\subsection{Model for Reconstructing the Coefficient $k$}
\label{subsec:model_k}

In this subsection, we propose a model for reconstructing the coefficient $k$.
The model is based on a MLP architecture, which leverages the universal approximation theorem~\cite{Cybenko89} as its theoretical foundation.
This theorem ensures that neural networks, including MLPs, can approximate any continuous function to a desired level of accuracy.

The neural network $\mathcal{N}$ takes the spatial variable $\xx$ as input and reconstructs the coefficient $k$ as $\tilde{k}$.
Specifically, the reconstruction is expressed as $\tilde{k}(\xx) = \mathcal{N}(\xx; \ttheta)$,
where $\ttheta$ represents the learnable parameters of the model, such as weights and biases.
The network consists of three hidden layers, each containing 50 hidden neurons.
To enable the model to capture complex natural signals and their derivatives, we use the sine function as the activation function in the hidden layers.
This periodic nonlinearity allows the model to represent fine details and high-frequency patterns more effectively than traditional activation functions~\cite{Sitzmann20}.

\subsection{Finite Element Method for the Forward Problem}
\label{subsec:forward}

In this subsection, we briefly describe the forward problem for the solution $p$ of \eqref{eq:Elliptic} by utilizing FEM.
Here, we employ the linear continuous Galerkin method, and it is implemented by FEniCS~\cite{Logg12}.
Given Dirichlet condition $g$ and coefficient $k$, \eqref{eq:Elliptic} can be reformulated into a weak form by integration and then approximated in a finite-dimensional space.

In particular, let $V_h$ be a finite-dimensional subspace of $H^1_0(\Omega)$ (to simplify the presentation, we assume $g=0$), and assume $p_h, v_h \in V_h$, with $p_h$ representing the approximate solution and $v_h$ being the test function used in the weak formulation.
Using the basis functions $\{\phi_i\}_{i = 1}^M$ of $V_h$, the weak form can be approximated as:
\begin{equation*}
\int_{\Omega} k\nabla p_h \cdot \nabla v_h \,{\rm d}\xx = \int_{\Omega} fv_h \,{\rm d}\xx, \ 
\forall v_h \in V_h.
\end{equation*}
Assuming $\displaystyle p_h = \sum_{i=1}^M P_i \phi_i$, 
where $P_i$ are constant coefficients, and substituting $v_h = \phi_j$, the resulting system becomes
\begin{equation*}
\sum_{i=1}^M P_i\int_{\Omega} k\nabla \phi_i \cdot \nabla \phi_j \,{\rm d}\xx = \int_{\Omega} f\phi_j \, {\rm d}\xx, \  \forall j=1,\cdots,M.
\end{equation*}
This can be written in matrix form as:
\begin{equation*}
\mathbf{A}\mathbf{P} = \mathbf{b},
\end{equation*}
where
\begin{equation*}
\mathbf{A}_{ij} = \int_{\Omega} k\nabla \phi_i \cdot \nabla \phi_j \,{\rm d}\xx, \quad \mathbf{P} = \begin{bmatrix} P_1 \\ \vdots \\ P_M \end{bmatrix}, \text{ and } \mathbf{b} = \begin{bmatrix} \int_{\Omega} f\phi_1 \,{\rm d}\xx \\ \vdots \\ \int_{\Omega} f\phi_M \,{\rm d}\xx \end{bmatrix}.
\end{equation*}
Here, $\mathbf{A}$ is known as the stiffness matrix, and $\mathbf{b}$ is the load vector.
Finally, solving $\mathbf{A} \mathbf{P} = \mathbf{b}$ yields $\mathbf{P}$, which determines the coefficients for the basis functions, allowing us to obtain the approximate solution $p_h$.

\subsection{Loss Functions}
\label{subsec:loss}

This subsection defines the loss function used to train the model introduced in Subsection~\ref{subsec:model_k}, and describes how to compute its gradient with respect to the model parameters.

The loss function based on \textbf{Problem} is defined as
\begin{equation}\label{eq:loss}
    \mathcal{L} (\mathcal D, \tilde{\mathcal D}) = \sum_{i = 1}^N \bigl\| h_{g_i} - \tilde{k} \partial_\nu \tilde{p}_{g_i}\bigr\|_{L^2(\omega)}^2,
\end{equation}
where $\mathcal D = \{(g_i, h_{g_i})_{i = 1, \dots, N} \mid g_i \in H^{1/2}(\partial\Omega), \, \operatorname{supp}g_i \subset \omega, h_{g_i} = k \partial_\nu p_{g_i}|_{\omega}\}$,  
$\tilde{\mathcal D} = \{(g_i, \tilde{k}\partial_\nu \tilde{p}_{g_i}|_{\omega})_{i = 1, \dots, N} \mid g_i \in H^{1/2}(\partial\Omega), \, \operatorname{supp}g_i \subset \omega\}$,  
and $\|\cdot\|_{L^2(\omega)}$ denotes the $L^2$ norm on $\omega$.
The loss function $\mathcal{L}$ is minimized by adjusting the model parameters $\ttheta$.
If the loss function reaches zero, then $\mathcal D = \tilde{\mathcal D}$.

To optimize the model, we need the gradient of the loss function $\mathcal{L}$ with respect to the parameters $\ttheta$.
However, since this gradient cannot be computed directly, we apply the chain rule:
\begin{equation}\label{eq:chain}
    \frac{\partial \mathcal{L}}{\partial \ttheta} = \frac{\partial \mathcal{L}}{\partial \tilde{k}_{\text{fem}}} \cdot \frac{\partial \tilde{k}}{\partial \ttheta}.
\end{equation}
Here, $\tilde{k}$ represents the neural network output as a PyTorch tensor and $\tilde{k}_{\text{fem}} = \tilde{k}.\texttt{detach()}.\texttt{numpy()}$ is a coefficient that shares the same values as $\tilde{k}$ but is converted to a NumPy array for use in the forward PDE solve within FEniCS.
The $\texttt{detach()}$ operation disconnects $\tilde{k}$ from PyTorch’s autograd graph, meaning that any gradient information associated with it is lost beyond this point.
As a result, the gradient $\partial \mathcal{L}/\partial \tilde{k}_{\text{fem}}$ is computed efficiently using the adjoint method in FEniCS, while $\partial \tilde{k}/\partial \ttheta$ is obtained via PyTorch’s automatic differentiation~\cite{Paszke19}.
Multiplying these two terms yields the desired gradient $\partial \mathcal{L}/\partial \ttheta$.

\section{Numerical Experiments}
\label{sec:num}
In this section, we present several numerical experiments to validate and demonstrate the effectiveness of the proposed framework introduced in Section~\ref{sec:algorithm}.
The following setup is common across all examples.
The neural network is trained for 100 epochs with a batch size of 1, using PyTorch~\cite{Paszke19} and the ADAM optimizer~\cite{Kingma14}, a variant of stochastic gradient descent.
All experiments are performed on a system equipped with an NVIDIA GeForce RTX 3090 GPU (24GB VRAM), an Intel Core i9-12900K CPU, and 62GB of RAM.
The software environment includes CUDA 11.8 and NVIDIA driver version 522.06.

\subsection{Example 1-1: Constant Coefficient}
In this example, we consider a constant coefficient $k$ to provide a baseline for evaluating model performance under simple conditions.
Specifically, we set $k = 1$.
The computational domain is the unit square $\Omega = (0,1)\times(0,1)$, discretized using a triangular mesh with 4,225 degrees of freedom for the linear Lagrange finite element space.
The Dirichlet boundary condition is defined as $g(x,y) = \cos(a\pi x)\cos(b\pi y)|_{\partial\Omega}$, where $a$ and $b$ are independently sampled from the uniform distribution over $[0,2]$.
This experiment uses a Cauchy dataset containing 2,048 samples.
The source function is given by $f(x,y) = -\nabla \cdot (k(x,y) \nabla (\cos(a\pi x)\cos(b\pi y)))$.
The learning rate is initialized at 0.001 and reduced by a factor of 0.25 every 20 epochs to ensure stable convergence and prevent overshooting in the later stages of training.

To evaluate the effectiveness of our method in reconstructing the coefficient $k$, we employ two measures: the Neumann data loss and the relative error of $\tilde{k}$.
The Neumann data loss, defined in \eqref{eq:loss} with~$\omega = \partial\Omega$ (fully observed setting), measures the difference between the predicted and observed Neumann data.
Meanwhile, the relative error is defined as
\begin{equation}\label{eq:relative_error}
    \|\tilde{k}\|_{\text{rel}} = \frac{\|k - \tilde{k}\|_{L^2(\overline{\Omega})}}{\|k\|_{L^2(\overline{\Omega})}},
\end{equation}
where $\|\cdot\|_{L^2(\overline{\Omega})}$ denotes the $L^2$ norm in $\overline{\Omega}$, and provides a normalized measure of the difference between $\tilde{k}$ and the exact coefficient $k$.

\begin{figure}[!ht]
    \centering

    \begin{subfigure}{0.55\textwidth}
        \centering
        \includegraphics[width=\textwidth]{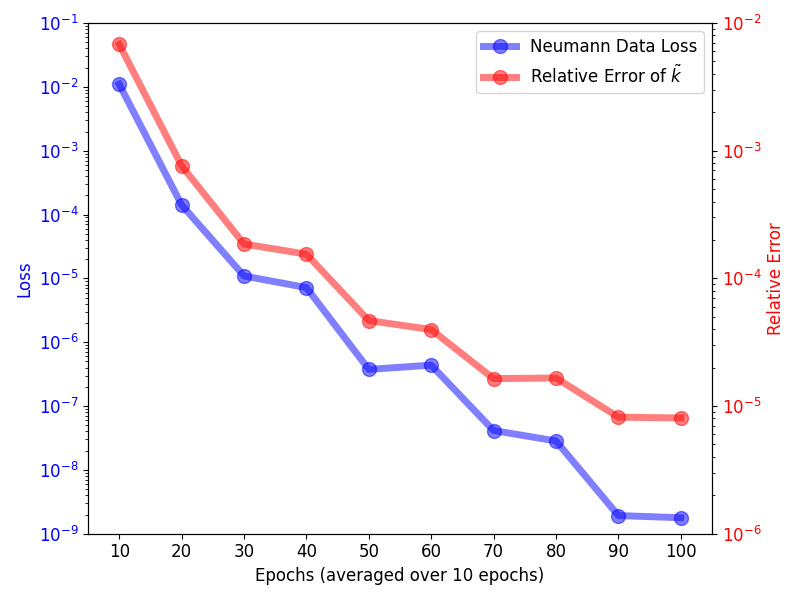}
        \caption{Neumann data loss and relative error of $\tilde{k}$}
    \end{subfigure}
    \hfill
    \begin{subfigure}{0.40\textwidth}
        \centering
        \includegraphics[width=\textwidth]{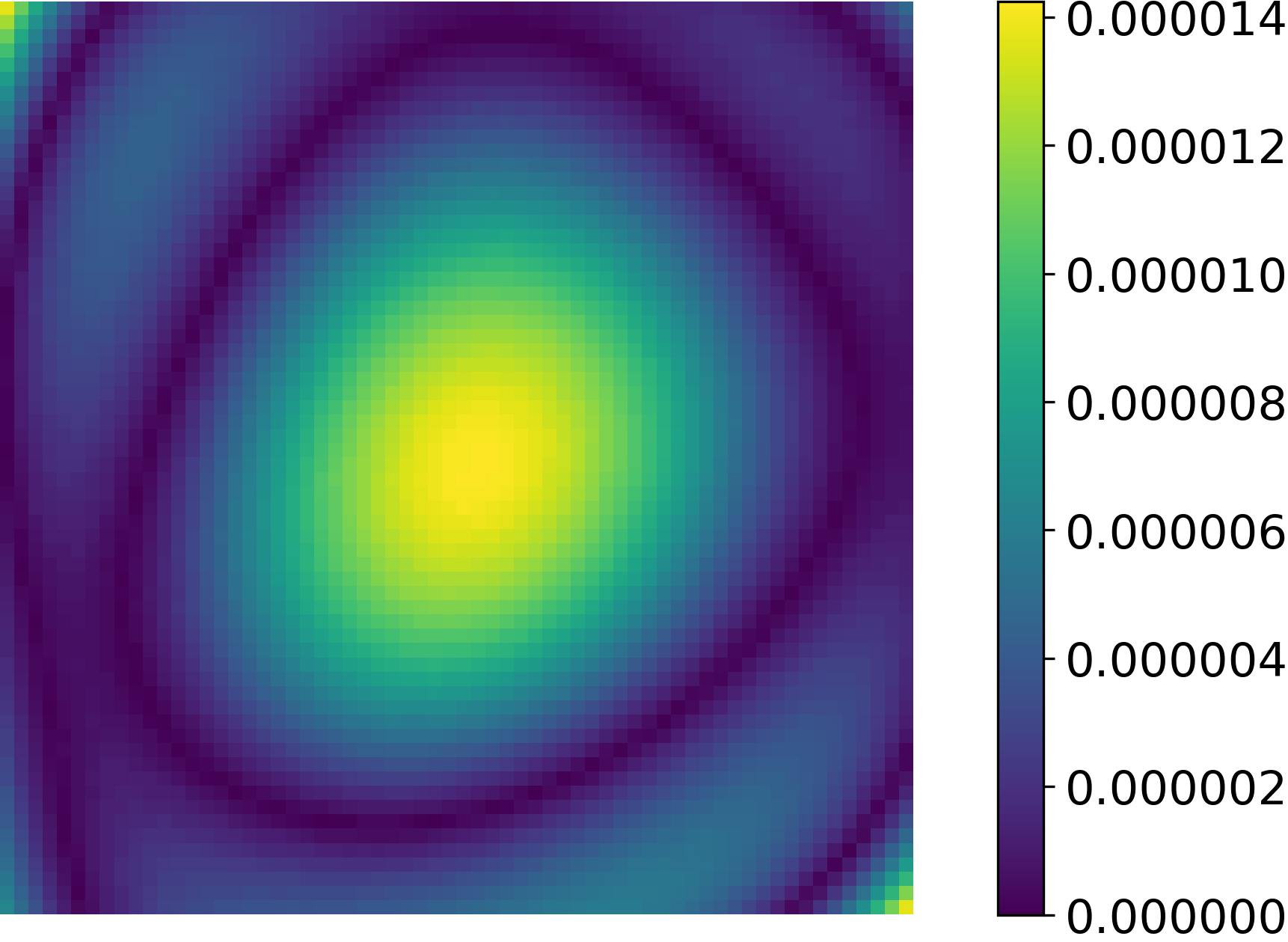}
        \caption{absolute error $|k - \tilde{k}|$.}
    \end{subfigure}

    \caption{Example 1-1. Training dynamics and reconstruction accuracy.
    (a) Neumann data loss (blue, left axis) and relative error of $\tilde{k}$ (red, right axis) over training epochs (log scale). Each point represents the averaged value over a 10-epoch interval.
    (b) Absolute error distribution between the exact coefficient $k$ and the reconstructed coefficient $\tilde{k}$.}
    \label{fig:constant_results}
\end{figure}
Figure~\ref{fig:constant_results} illustrates both the training dynamics and the reconstruction accuracy.
Figure~\ref{fig:constant_results}(a) shows the Neumann data loss (blue, left axis) and the relative error of $\tilde{k}$ (red, right axis), both plotted on a logarithmic scale.
Each point corresponds to the average over a 10-epoch interval, which smooths out short-term fluctuations and highlights the overall decreasing trend.
The simultaneous reduction of both metrics demonstrates that minimizing the Neumann data loss effectively drives $\tilde{k}$ toward the exact coefficient $k$, supporting the uniqueness claim in \textbf{Problem}.
Figure~\ref{fig:constant_results}(b) visualizes the absolute error $|k - \tilde{k}|$ between the exact coefficient $k$ and the reconstructed coefficient $\tilde{k}$.
The results confirm that the proposed method provides a highly accurate approximation of $k$, with only minor deviations.

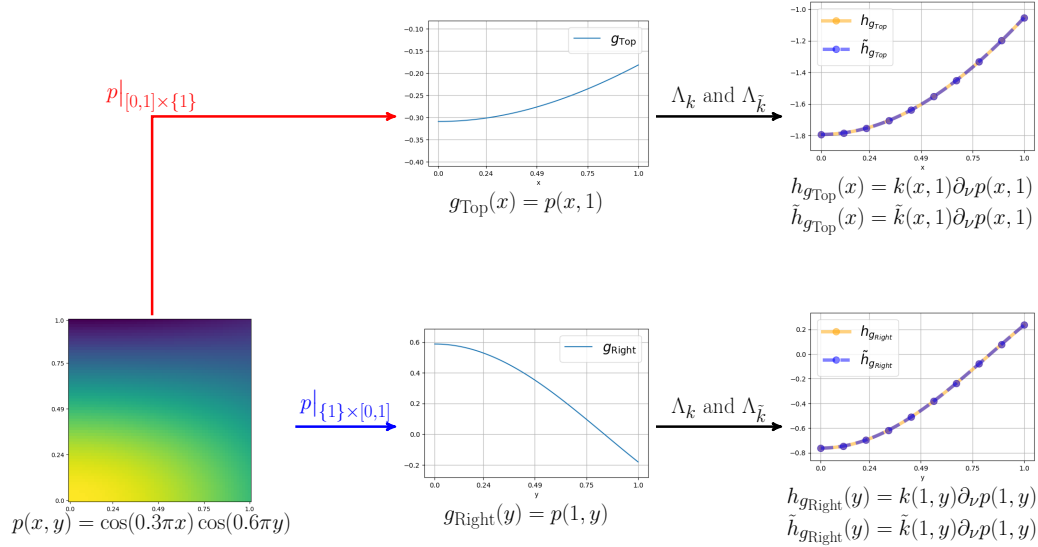
\begin{figure}[!ht]
    \centering
    \resizebox{\textwidth}{!}{

\begin{tikzpicture}[>=stealth, node distance=3cm, font=\small]
    \tikzset{line/.style={draw,->,>=stealth'}}
      \node[draw=none, rectangle, minimum width=5cm, minimum height=2cm, align=center] (p) 
      {%
        \includegraphics[width=8cm]{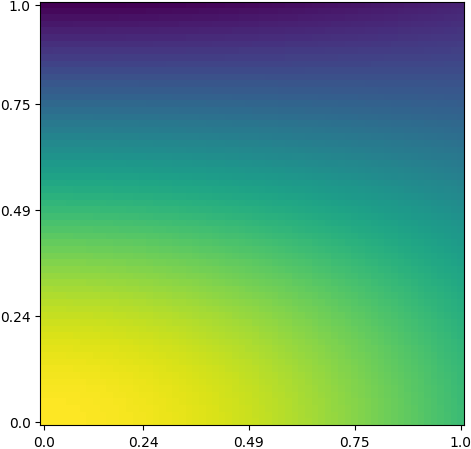}\\
        \huge $p(x,y) = \cos(0.3 \pi x) \cos(0.6 \pi y)$
      };

      \node[draw=none, rectangle, minimum width=10cm, minimum height=2cm, align=center,
            right=4cm of p] (p_Right)
      {%
        \includegraphics[width=10cm]{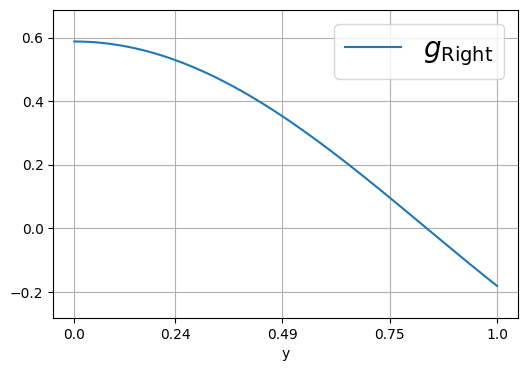}\\
        \huge $g_{\text{Right}}(y) = p(1,y)$
      };
  
      \node[draw=none, rectangle, minimum width=10cm, minimum height=2cm, align=center,
      above=4cm of p_Right] (p_Top)
      {%
      \includegraphics[width=10cm]{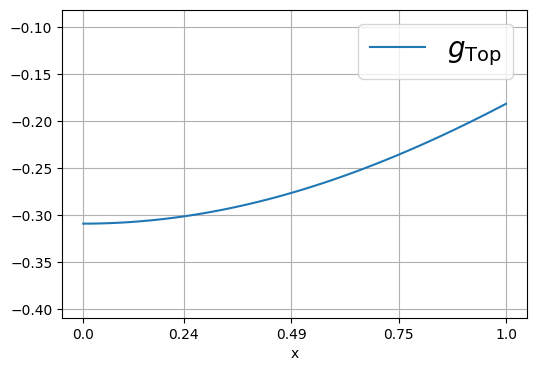}\\
          \huge $g_{\text{Top}}(x) = p(x,1)$
      };
  
      \node[draw=none, rectangle, minimum width=3.0cm, minimum height=2cm, align=center,
            right=5.0cm of p_Top] (Neumann_Top)
    {%
      \includegraphics[width=10cm]{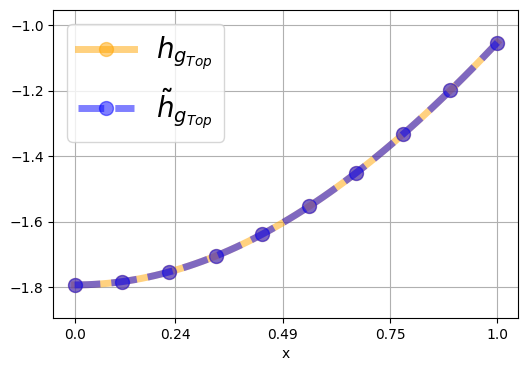}\\
      \huge $h_{g_{\text{Top}}}(x) = k(x,1)\partial_\nu p(x,1)$\\
      \huge $\tilde{h}_{g_{\text{Top}}}(x) = \tilde{k}(x,1)\partial_\nu p(x,1)$
    };
  
    \node[draw=none, rectangle, minimum width=3.0cm, minimum height=2cm, align=center,
          right=5.0cm of p_Right] (Neumann_Right)
    {%
      \includegraphics[width=10cm]{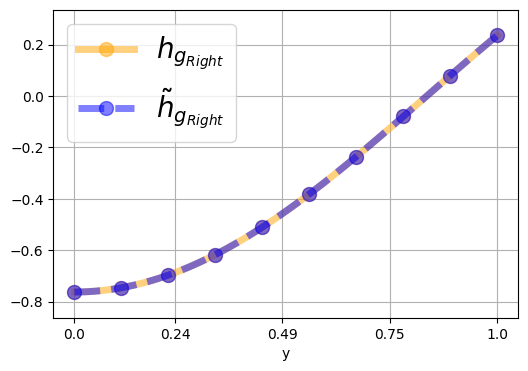}\\
      \huge $h_{g_{\text{Right}}}(y) = k(1,y)\partial_\nu p(1,y)$\\
      \huge $\tilde{h}_{g_{\text{Right}}}(y) = \tilde{k}(1,y)\partial_\nu p(1,y)$
    
    };

    \path[line, line width=1mm, red] (p.north)  |- (p_Top.west) node[midway, above] {\Huge $p \big|_{[0,1] \times \{1\}}$};
    \path[line, line width=1mm, blue] (p.east)  -- (p_Right.west) node[midway, above] {\Huge $p \big|_{\{1\} \times [0,1]}$};
    \path[line, line width=1mm] (p_Top.east) -- (Neumann_Top.west) node[midway, above] {\Huge $\Lambda_k$ and $\Lambda_{\tilde{k}}$};
    \path[line, line width=1mm] (p_Right.east) -- (Neumann_Right.west) node[midway, above] {\Huge $\Lambda_k$ and $\Lambda_{\tilde{k}}$};
  
  \end{tikzpicture}%
  }
  \caption{Example 1-1. Illustration of the Dirichlet-to-Neumann (DtN) map. The solution $p(x,y)$ is shown and the Dirichlet boundary values $p(x,1)$ and $p(1,y)$ mapped to their corresponding Neumann data via the DtN map are presented.}
\label{fig:illustrate} 

\end{figure}

To provide a more intuitive understanding of how the boundary conditions operate in practice, we now examine the operation of the DtN map as illustrated in Figure~\ref{fig:illustrate}.
The heat map on the left represents the solution $p(x,y)=\cos(0.3\pi x)\cos(0.6\pi y)$  over a rectangular domain.
From this solution, Dirichlet data are obtained along the top and right boundaries, given by $g_{\text{Top}}(x) = p(x,1)$ and $g_{\text{Right}}(y) = p(1,y)$.
The DtN map converts these Dirichlet data into the corresponding Neumann data, $h_{g_{\text{Top}}}(x)$ and $h_{g_{\text{Right}}}(y)$, as shown in the adjacent plots.

\subsection{Example 1-2: Constant Coefficient with Partially Observed Data}
In this example, we address the inverse problem under the partially observed data.
In many real-world scenarios, measurements may be unavailable in certain regions due to sensor limitations or other constraints.

We consider the case of a constant coefficient $k=1$. All other configurations, including the computational domain, degrees of freedom, the number of data samples, and the learning rate, are the same as in Example 1-1.
To simulate partial observations, we assume that measurements are unavailable in certain parts of the boundary, covering approximately 40\% of $\partial \Omega$.
Specifically, we define the region $\Gamma_0$, where measurements are unavailable, as follows:
measurements are unavailable for $y \in [\frac{2}{5}, \frac{4}{5}]$ on the left boundary ($x = 0$), for $y \in [\frac{1}{5}, \frac{3}{5}]$ on the right boundary ($x = 1$), for $x \in [\frac{1.5}{5}, \frac{3.5}{5}]$ on the bottom boundary ($y = 0$), and for $x \in [\frac{2.5}{5}, \frac{4.5}{5}]$ on the top boundary ($y = 1$).
In these regions, both the Dirichlet boundary condition and the source function are set to zero; that is, $g = 0$ and $f = 0$ on $\Gamma_0$.
The remaining boundary, denoted by $\Gamma = \partial \Omega \setminus \Gamma_0$, is where measurements are available.
On $\Gamma$, the Dirichlet boundary condition is defined as $g(x,y) = \cos(a\pi x)\cos(b\pi y)$, where $a$ and $b$ are independently sampled from a uniform distribution over $[0,2]$.
Similarly, the source function is defined on $\Gamma$ as  $f(x,y) = -\nabla \cdot (k(x,y) \nabla (\cos(a\pi x)\cos(b\pi y)))$.
The Neumann data is measured only on $\Gamma$, and the Neumann data loss function~\eqref{eq:loss} is defined with $\omega = \Gamma$.
\begin{figure}[!ht]
    \centering
    \includegraphics[width=0.4\textwidth]{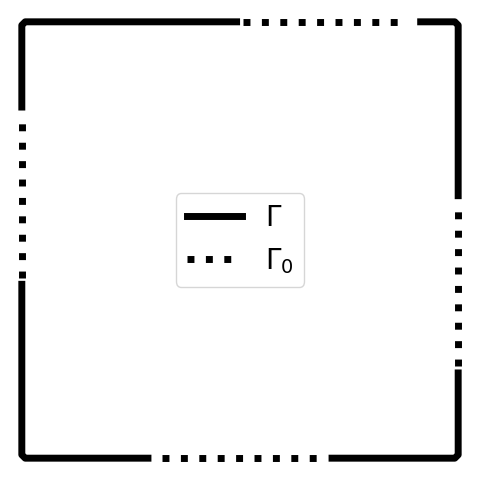}
    \includegraphics[width=0.48\textwidth]{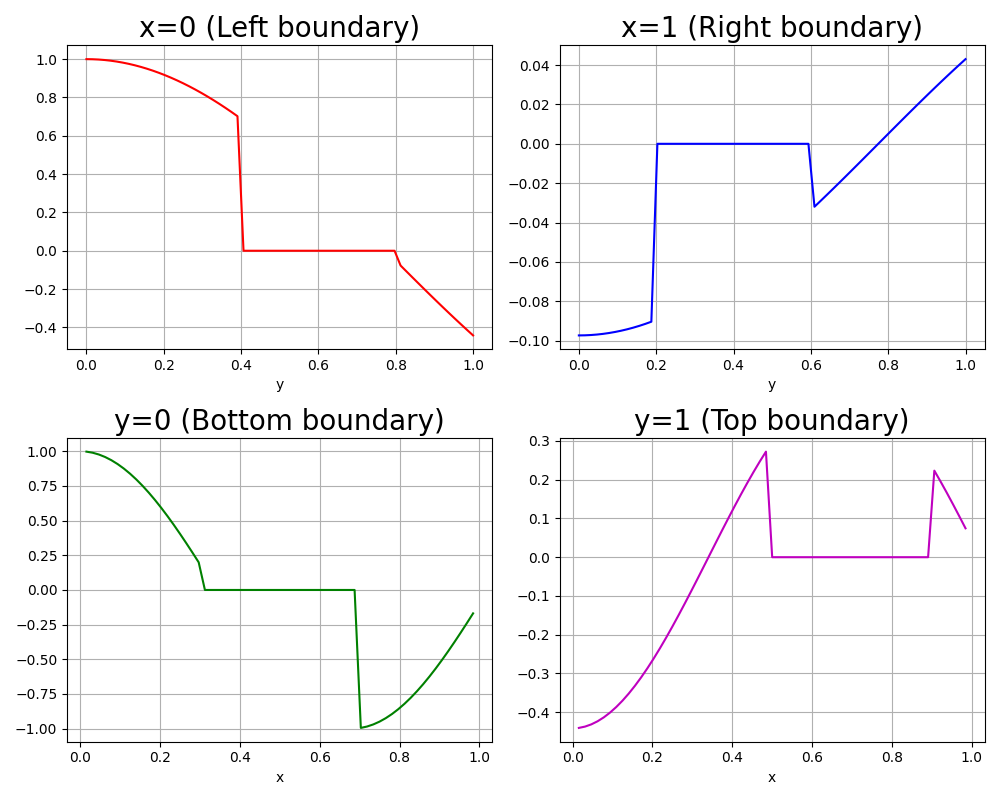}
    \caption{Example 1-2. (Left) Boundary regions and (Right) Dirichlet boundary condition $g$. The values for which  $g=0$ indicate that no data is available.}
    \label{fig:boundary}
\end{figure}
Figure~\ref{fig:boundary} illustrates the boundary region (left) and the corresponding Dirichlet boundary condition $g$ (right).

\begin{figure}[!ht]
    \centering

    \begin{subfigure}{0.55\textwidth}
        \centering
        \includegraphics[width=\textwidth]{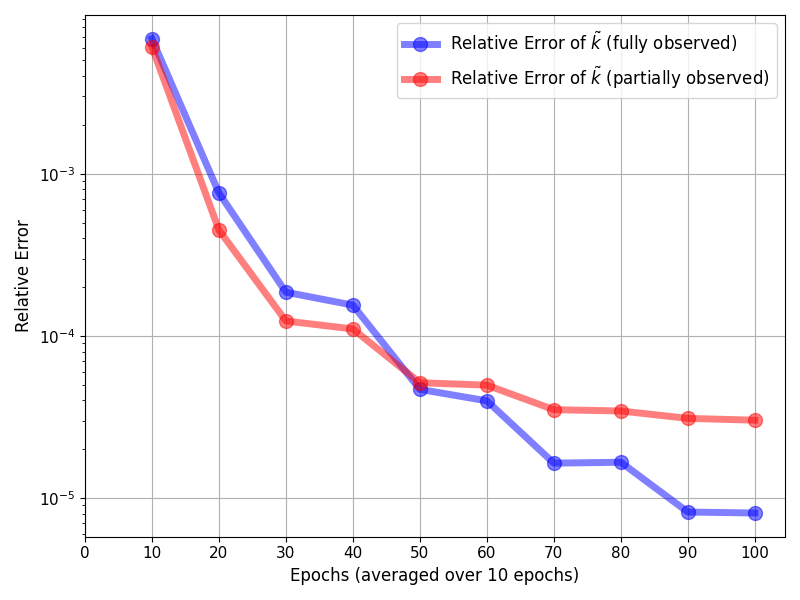}
        \caption{Relative error of $\tilde{k}$ over epochs (log scale) for fully and partially observed settings.}
    \end{subfigure}
    \hfill
    \begin{subfigure}{0.40\textwidth}
        \centering
        \includegraphics[width=\textwidth]{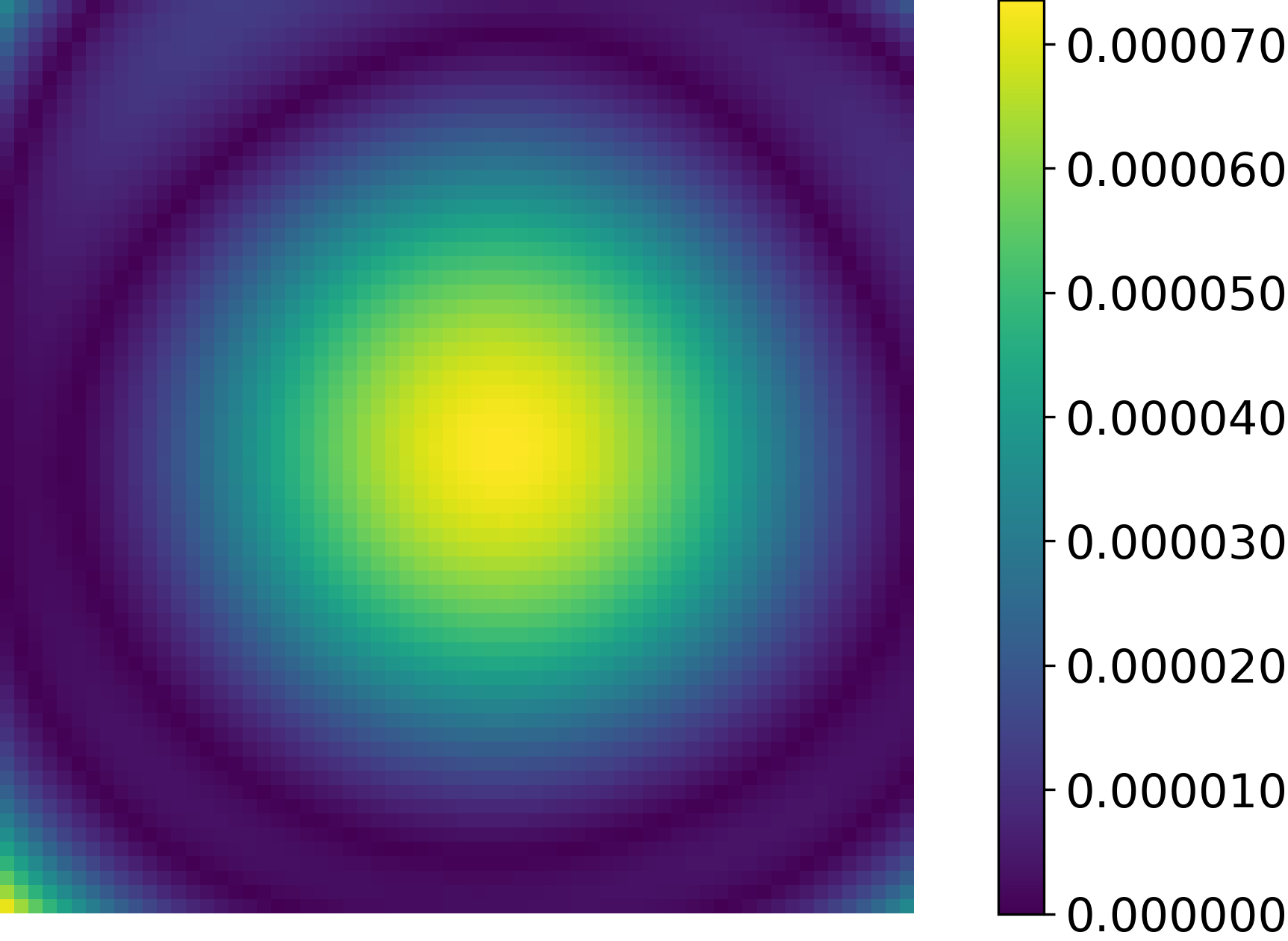}
        \caption{absolute error $|k - \tilde{k}|$.}
    \end{subfigure}

    \caption{Example 1-2. Training results and reconstruction accuracy for both fully and partially observed settings.
    (a) Relative error curves of $\tilde{k}$ under different observation settings.
    (b) Absolute error $|k - \tilde{k}|$ in the partially observed case.}
    \label{fig:constant_disconti_combined}
\end{figure}
Figure~\ref{fig:constant_disconti_combined} illustrates both the training results and the reconstruction accuracy.
Figure~\ref{fig:constant_disconti_combined}(a) shows the relative error of $\tilde{k}$ on a logarithmic scale, comparing the fully observed ($\omega = \partial\Omega$) and partially observed ($\omega = \Gamma$) settings.
Although the relative error in the partially observed setting is slightly higher than that in the fully observed case, it steadily decreases to a low value, indicating effective reconstruction.
Figure~\ref{fig:constant_disconti_combined}(b) visualizes the absolute error $|k - \tilde{k}|$ between the exact coefficient $k$ and the reconstructed coefficient $\tilde{k}$ in the partially observed setting.
These results indicate that the proposed method can accurately reconstruct the coefficient even with partial boundary measurements, demonstrating its robustness in handling incomplete data.

\subsection{Example 2-1: Spatially Varying Coefficient}
In this example, we test the algorithm with a spatially varying function $k$.
Specifically, we set
\begin{equation*}
k(x,y) = 0.9 \sin\left(\frac{\pi}{3}(x + y + 0.1)\right).
\end{equation*}
 All other configurations, including the computational domain, degrees of freedom, Dirichlet boundary condition, number of data samples, and learning rate, are the same as in Example 1-1.
The source function is defined as 
$f(x,y) = -\nabla \cdot (k(x,y) \nabla (\cos(a\pi x)\cos(b\pi y)))$
by computing the divergence of the flux associated with the given function.

\begin{figure}[!ht]
    \centering
    \begin{subfigure}{0.3\textwidth}
        \centering
        \includegraphics[width=\textwidth]{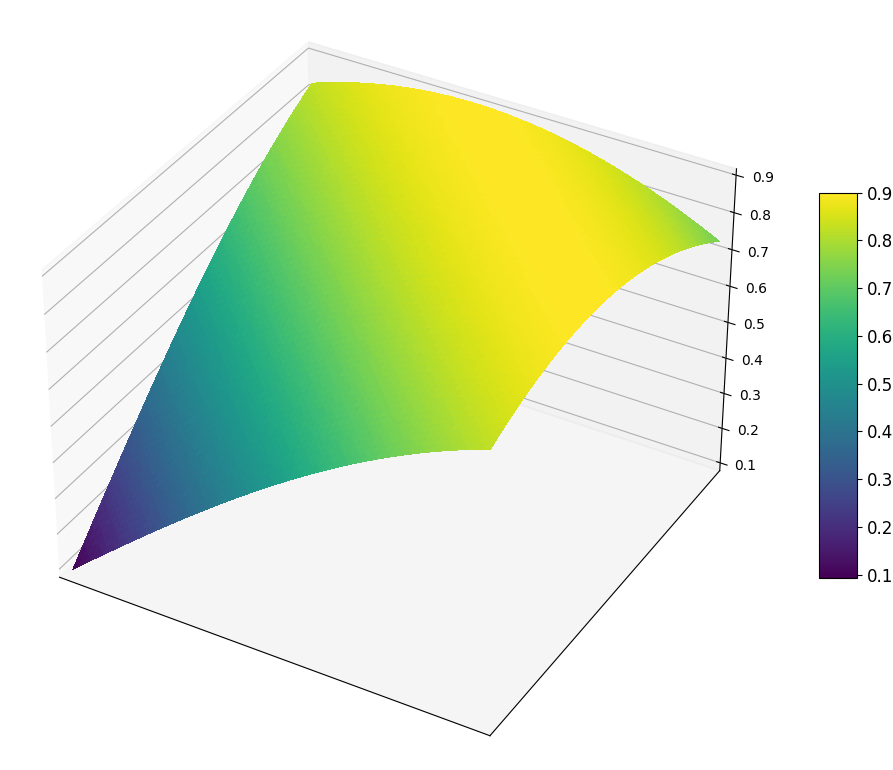}
        \caption{exact $k$}
    \end{subfigure}
    \hfill
    \begin{subfigure}{0.3\textwidth}
        \centering
        \includegraphics[width=\textwidth]{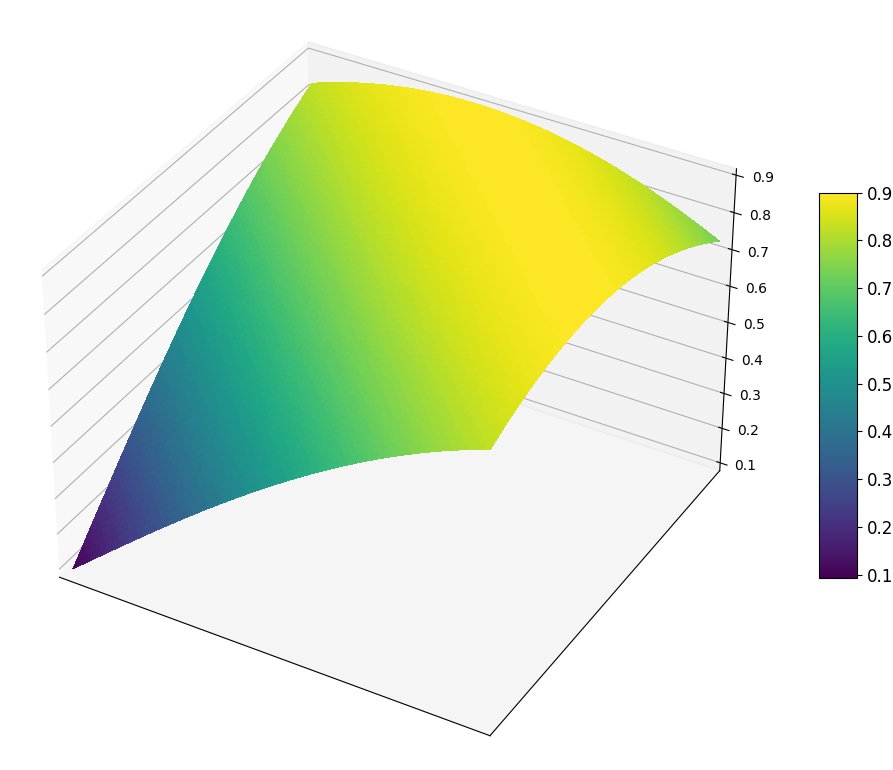}
        \caption{reconstructed $\tilde{k}$}
    \end{subfigure}
    \hfill
    \begin{subfigure}{0.33\textwidth}
        \centering
        \includegraphics[width=\textwidth]{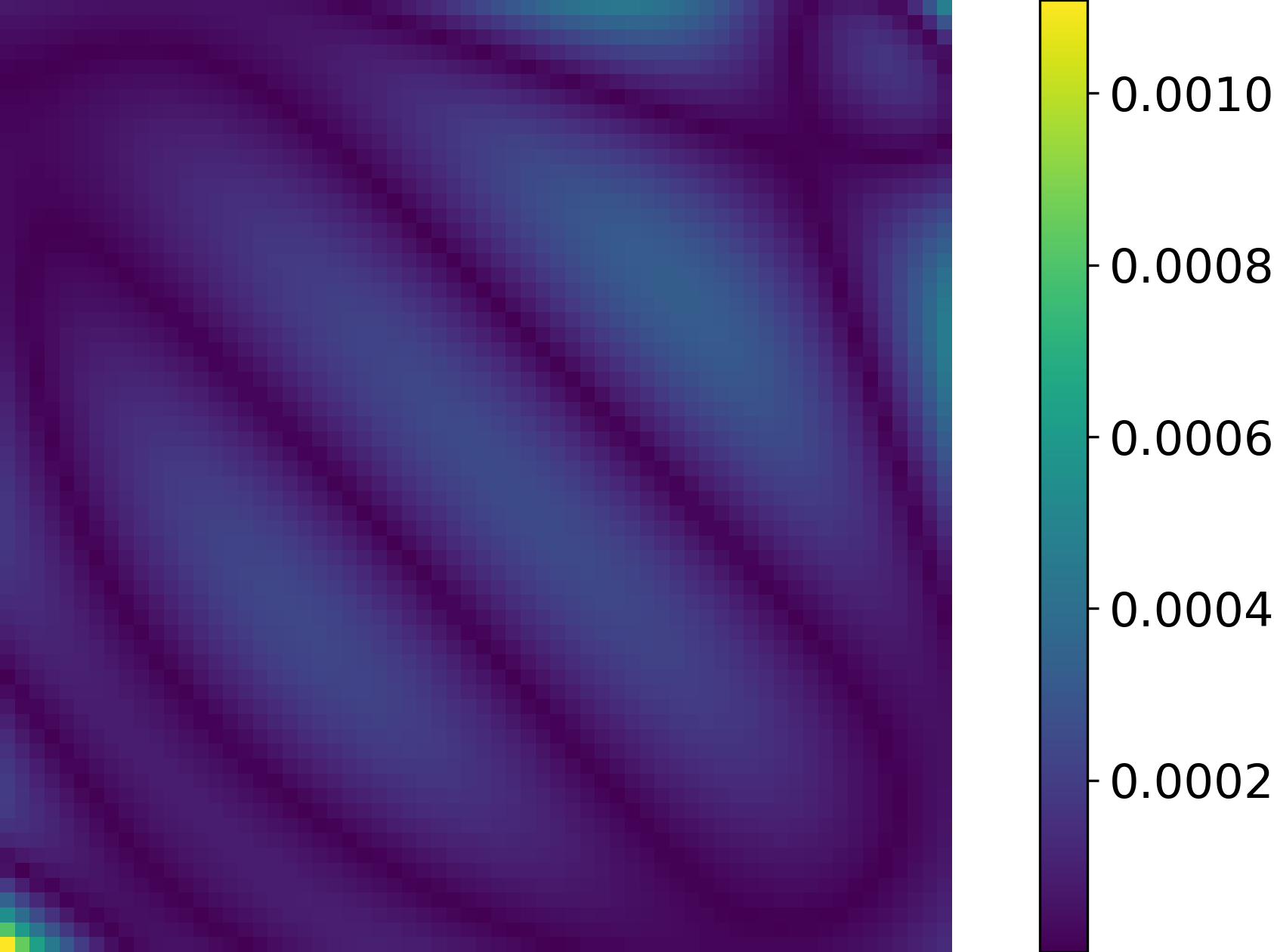}
        \vspace{-0.25cm} 
        \caption{absolute error $|k - \tilde{k}|$}
    \end{subfigure}

    \caption{Example 2-1. Visualizations of the exact coefficient $k$ and the reconstructed coefficient $\tilde{k}$ using 3D surface plots, along with the absolute error $|k - \tilde{k}|$ shown as a 2D heat map.}
    \label{fig:function_k}
\end{figure}
Figure \ref{fig:function_k} illustrates the exact coefficient $k$, the reconstructed coefficient $\tilde{k}$, and their absolute error, providing a comparative view of spatial structures.
The reconstructed coefficient $\tilde{k}$ closely matches the exact $k$, effectively capturing its spatial variations, with the absolute error remaining at a low level.

\begin{figure}[!ht]
    \centering
    \begin{subfigure}{0.45\textwidth}
        \centering
        \includegraphics[width=\textwidth]{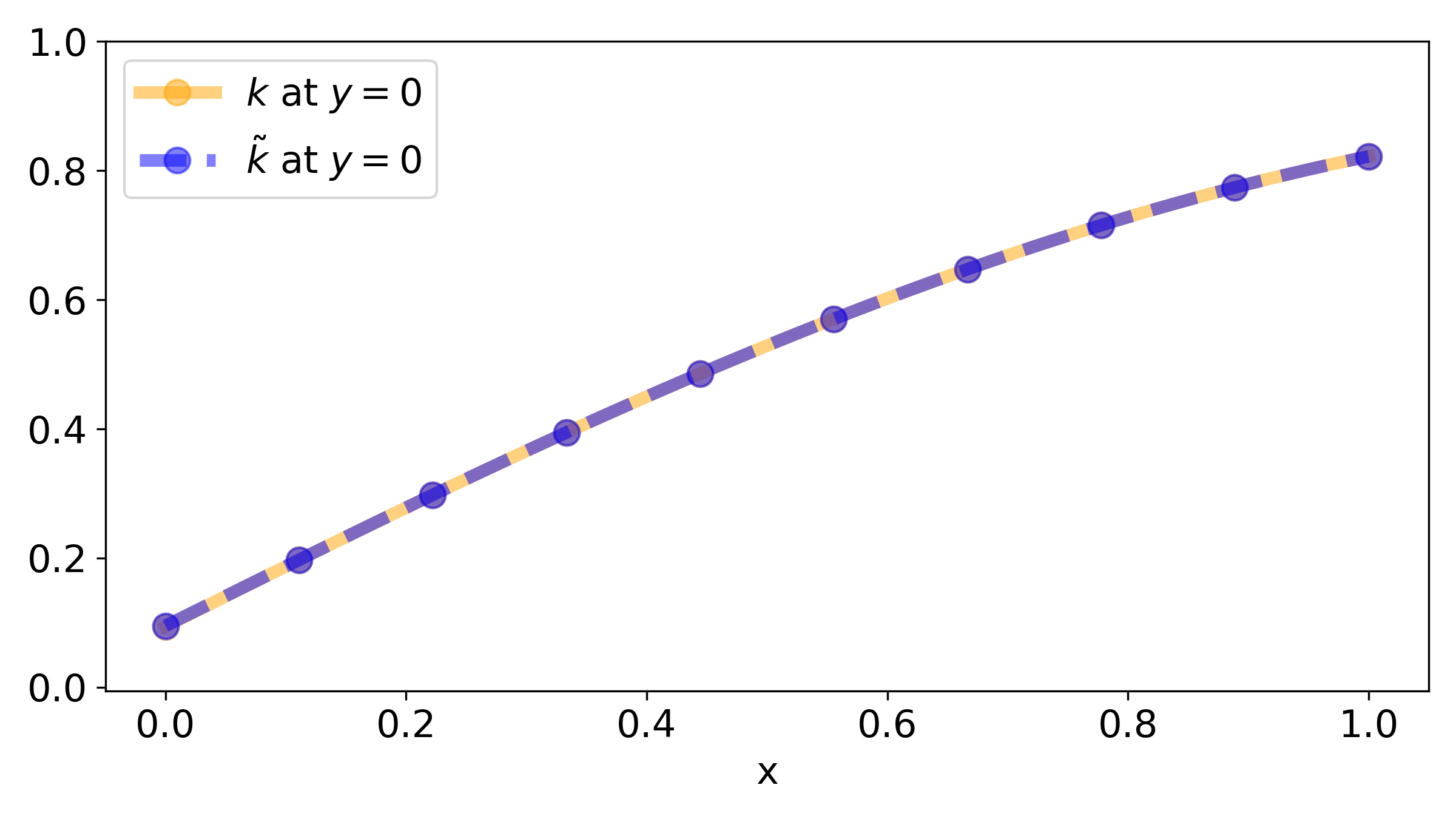}
        \caption{$y = 0$}
    \end{subfigure}
    \hfill
    \begin{subfigure}{0.45\textwidth}
        \centering
        \includegraphics[width=\textwidth]{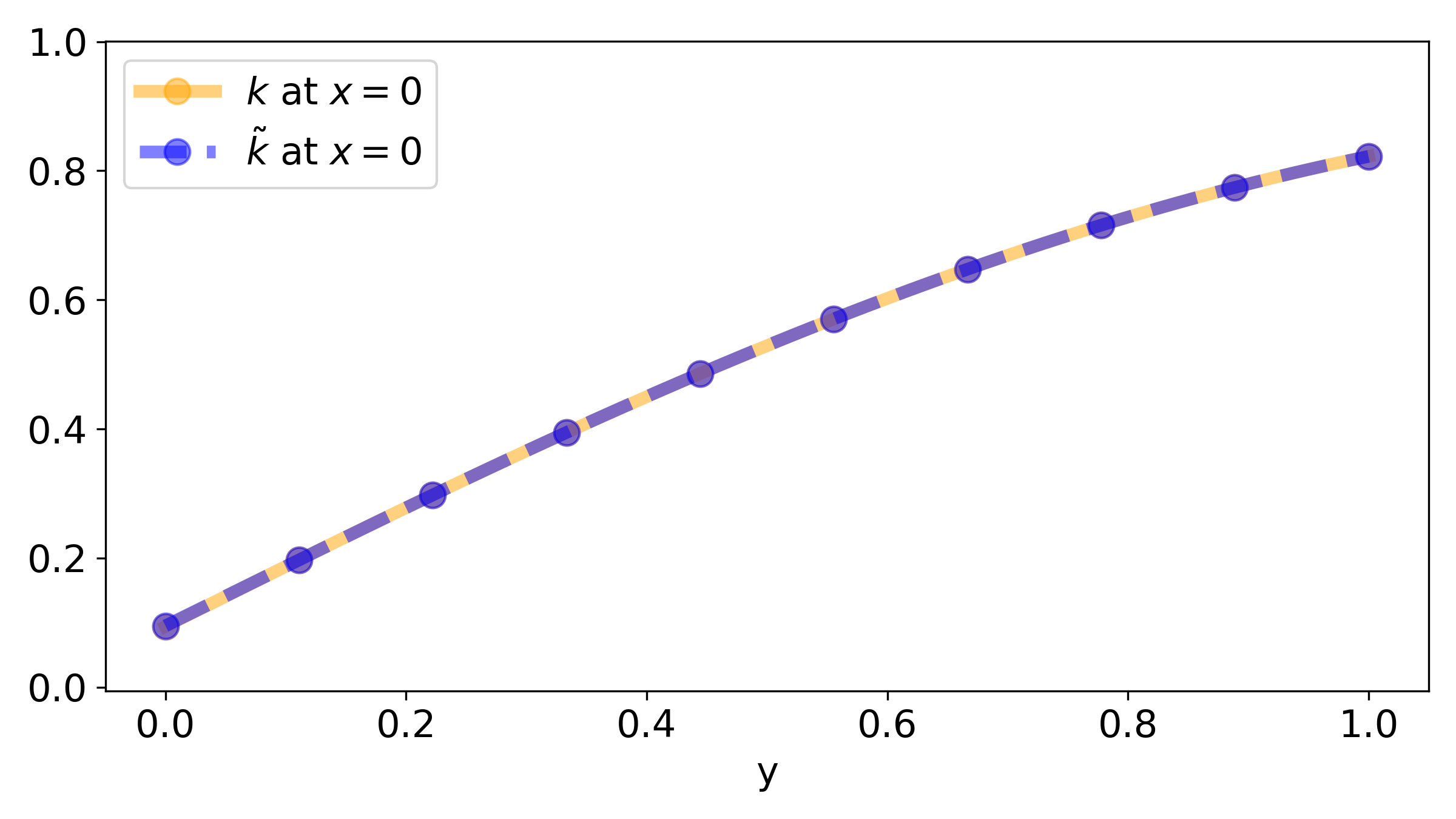}
        \caption{$x = 0$}
    \end{subfigure}
    \hfill
    \begin{subfigure}{0.45\textwidth}
        \centering
        \includegraphics[width=\textwidth]{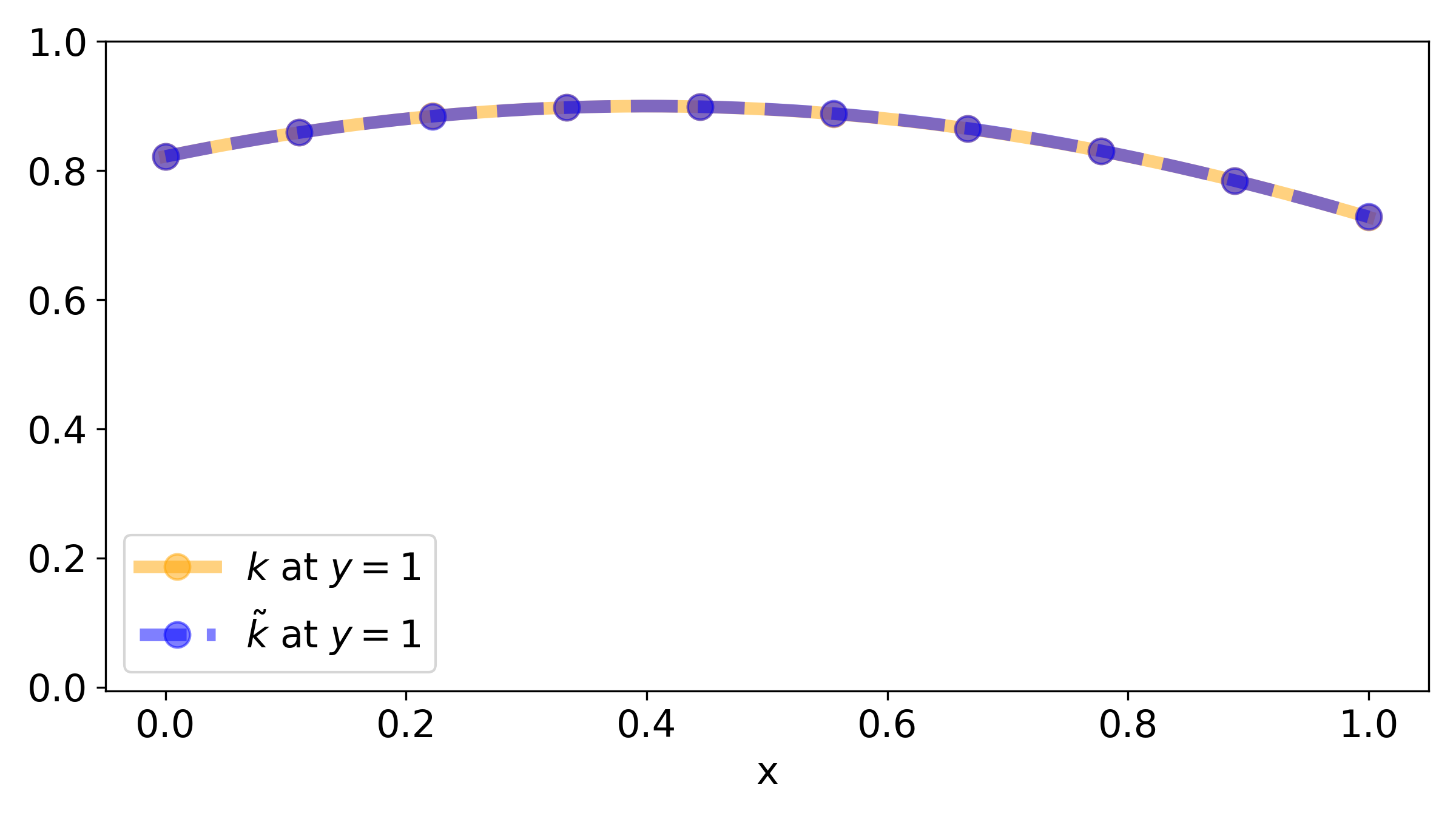}
        \caption{$y = 1$}
    \end{subfigure}
    \hfill
    \begin{subfigure}{0.45\textwidth}
        \centering
        \includegraphics[width=\textwidth]{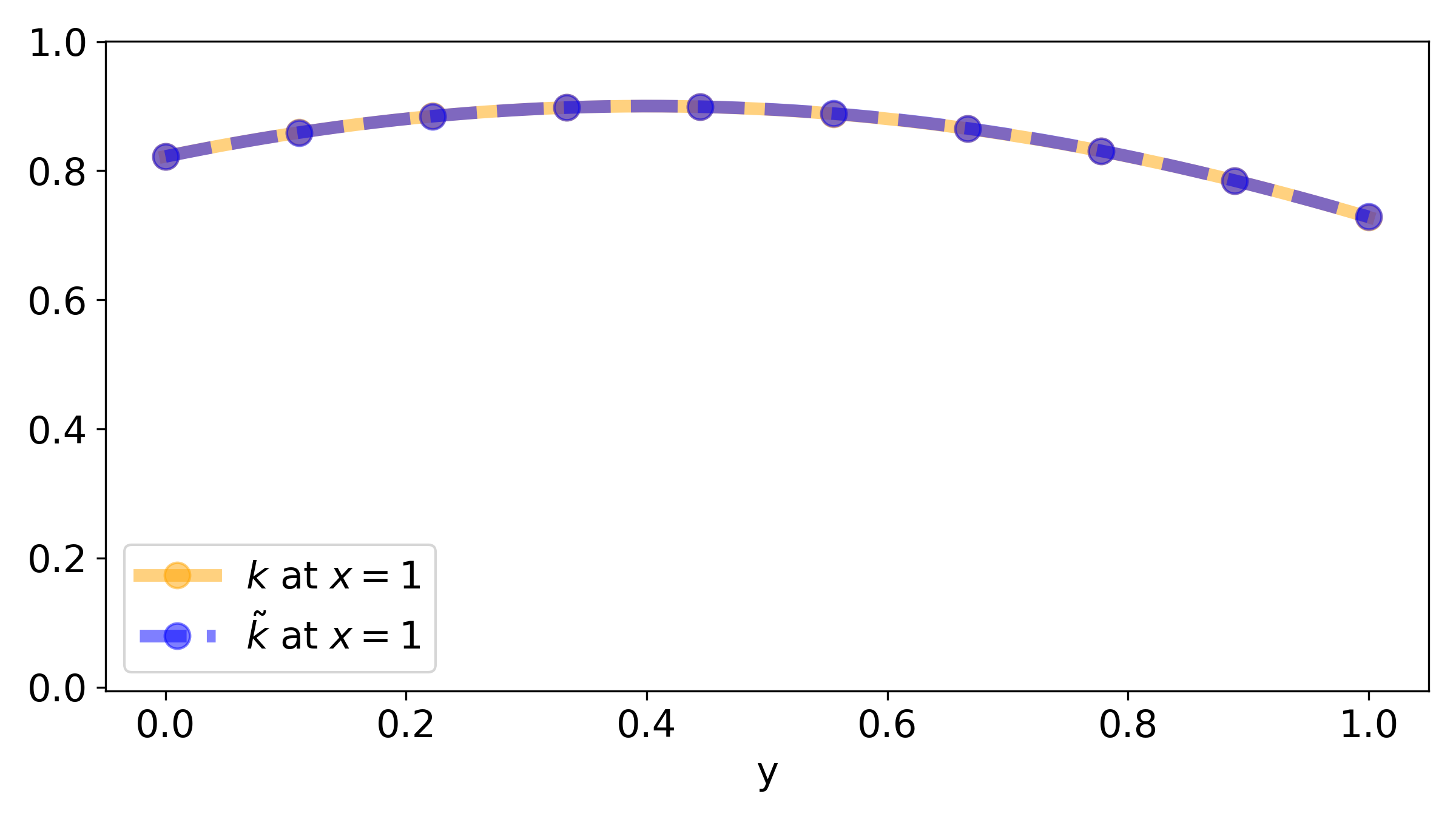}
        \caption{$x = 1$}
    \end{subfigure}
    \hfill
    \begin{subfigure}{0.45\textwidth}
        \centering
        \includegraphics[width=\textwidth]{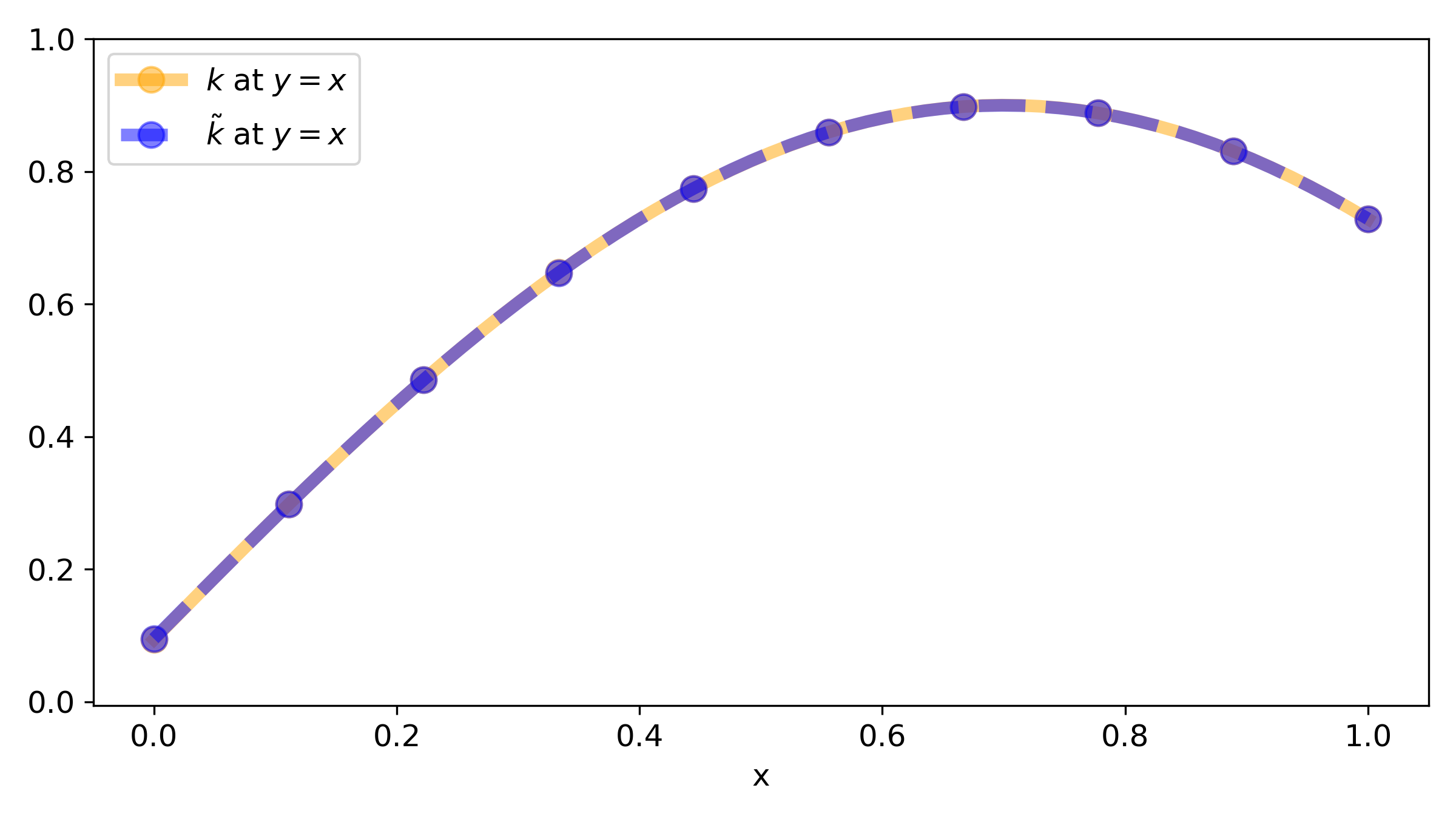}
        \caption{$y = x$}
    \end{subfigure}
    \hfill
    \begin{subfigure}{0.45\textwidth}
        \centering
        \includegraphics[width=\textwidth]{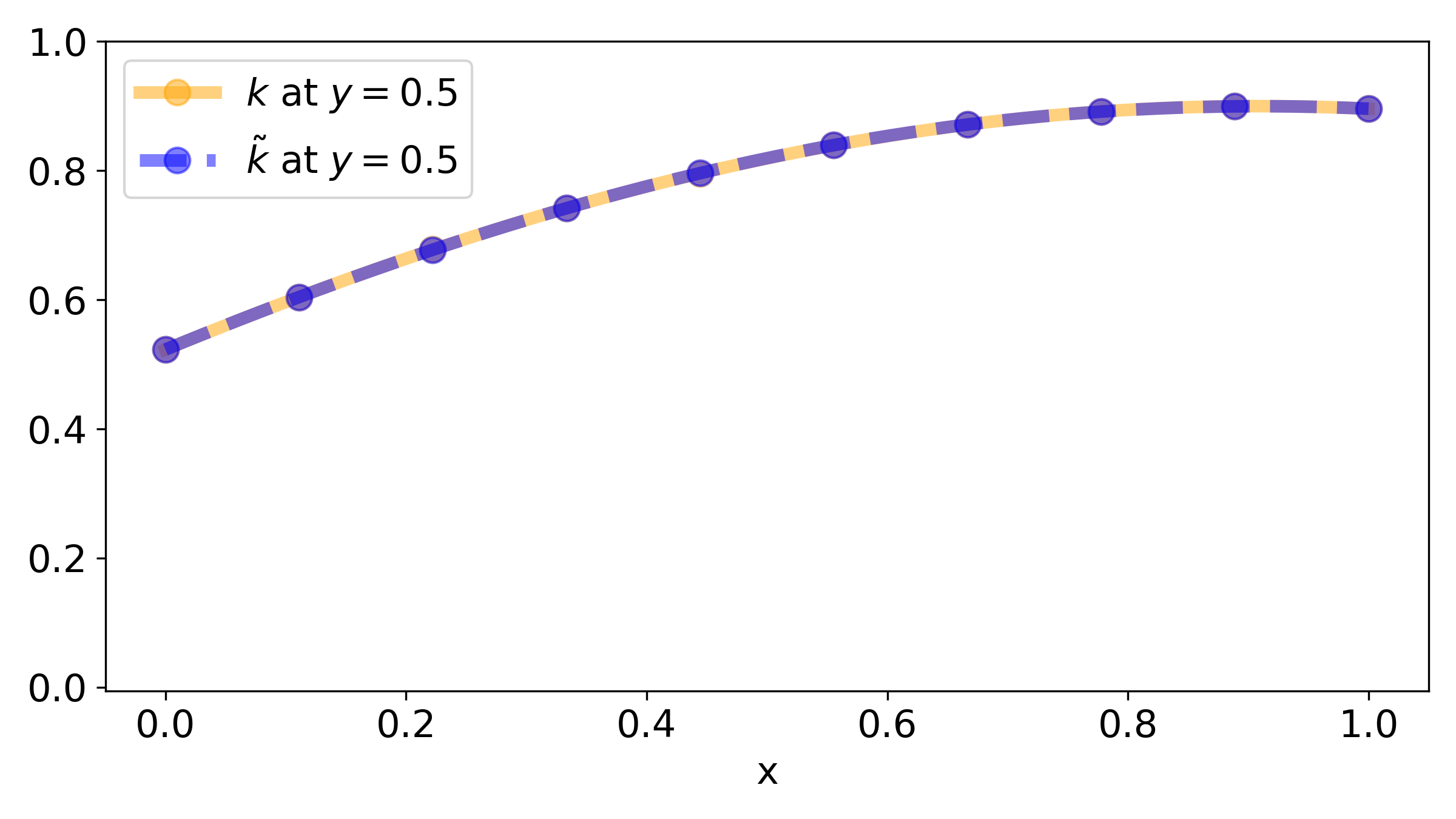}
        \caption{$y = 0.5$}
    \end{subfigure}

    \caption{Example 2-1. Slices of the exact coefficient $k$ and the reconstructed coefficient $\tilde{k}$ along selected lines in the domain.}
    \label{fig:function_k_slice}
\end{figure}

To further analyze the reconstruction accuracy, 
Figure \ref{fig:function_k_slice} presents slices of the coefficient $\tilde{k}$ and $k$ at the domain boundaries ($y=0$, $x=0$, $y=1$, and $x=1$) as well as interior slices along $y=x$ and $y=0.5$.
Figure~\ref{fig:function_k_slice} shows the reconstruction results by comparing both the boundary and interior regions.
Over the entire domain, the reconstructed coefficient $\tilde{k}$ closely follows the behavior of $k$, showing strong agreement both at the boundaries and the interior.
This consistency highlights the robustness of the method in accurately capturing spatial variations.

\subsection{Example 2-2: Spatially Varying Coefficient with Partially Observed Data}
In this example, we consider the case with partially observed data.
The coefficient $k$ is defined as $k(x,y) = 0.9 \sin(\frac{\pi}{3}(x + y + 0.1))$.
All other configurations, including the computational domain, degrees of freedom, Dirichlet boundary condition, number of data samples, and learning rate, are the same as in Example 1-2.
The source function is defined as 
$f(x,y) = -\nabla \cdot (k(x,y) \nabla (\cos(a\pi x)\cos(b\pi y)))$
by computing the divergence of the flux associated with the given function.
This function is imposed on $\Gamma$, while $f$ is set to zero on $\Gamma_0$.

\begin{figure}[ht!]
    \centering
    \includegraphics[width=0.6\textwidth]{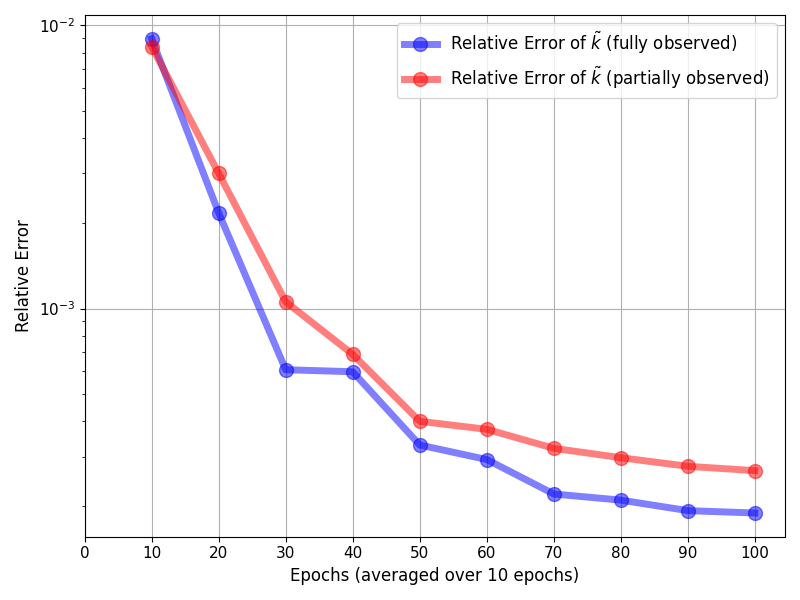}
    \caption{Example 2-2. Relative error of $\tilde{k}$ over epochs (log scale) for both fully and partially observed settings.  
    Each point represents the average over a 10-epoch interval.}
    \label{fig:function_disconti_combined}
\end{figure}
Figure~\ref{fig:function_disconti_combined} shows the relative error of $\tilde{k}$ on a logarithmic scale, comparing the fully observed and partially observed settings.
As in the constant coefficient case, the relative error under the partially observed setting is slightly higher, but it still decreases steadily, indicating good convergence behavior.
\begin{figure}[!ht]
    \centering

    \begin{subfigure}{0.3\textwidth}
        \centering
        \includegraphics[width=\textwidth]{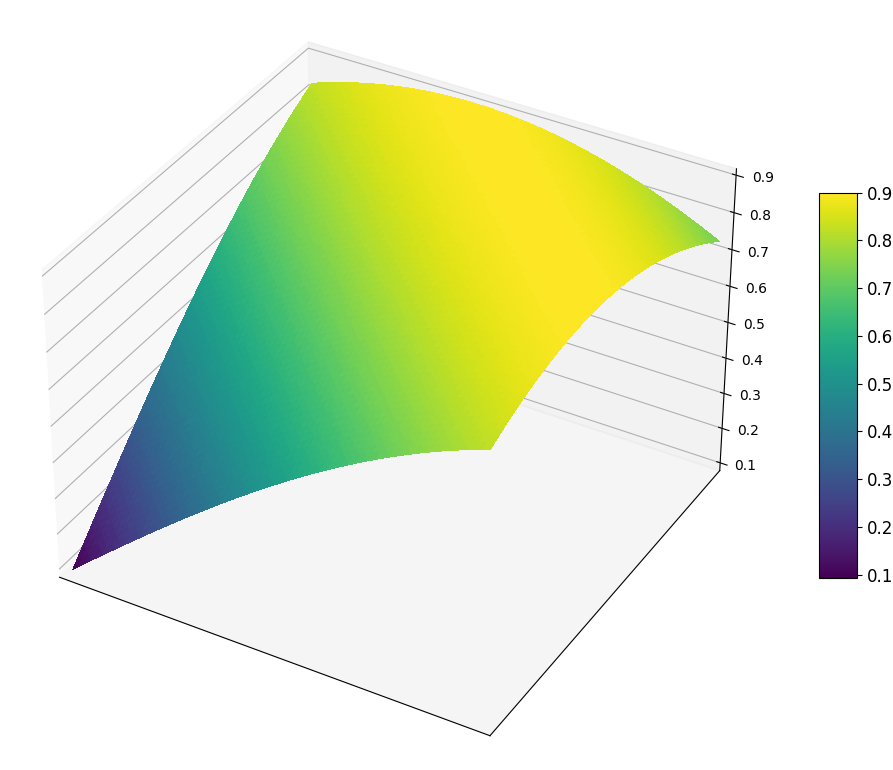}
        \caption{exact $k$}
    \end{subfigure}
    \hfill
    \begin{subfigure}{0.3\textwidth}
        \centering
        \includegraphics[width=\textwidth]{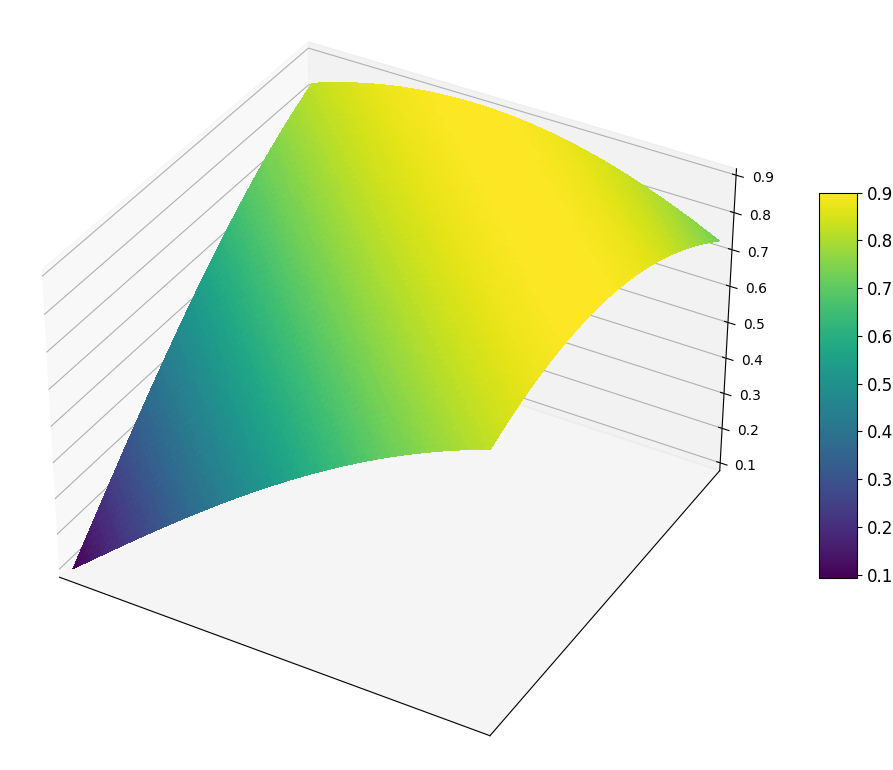}
        \caption{reconstructed $\tilde{k}$}
    \end{subfigure}
    \hfill
    \begin{subfigure}{0.33\textwidth}
        \centering
        \includegraphics[width=\textwidth]{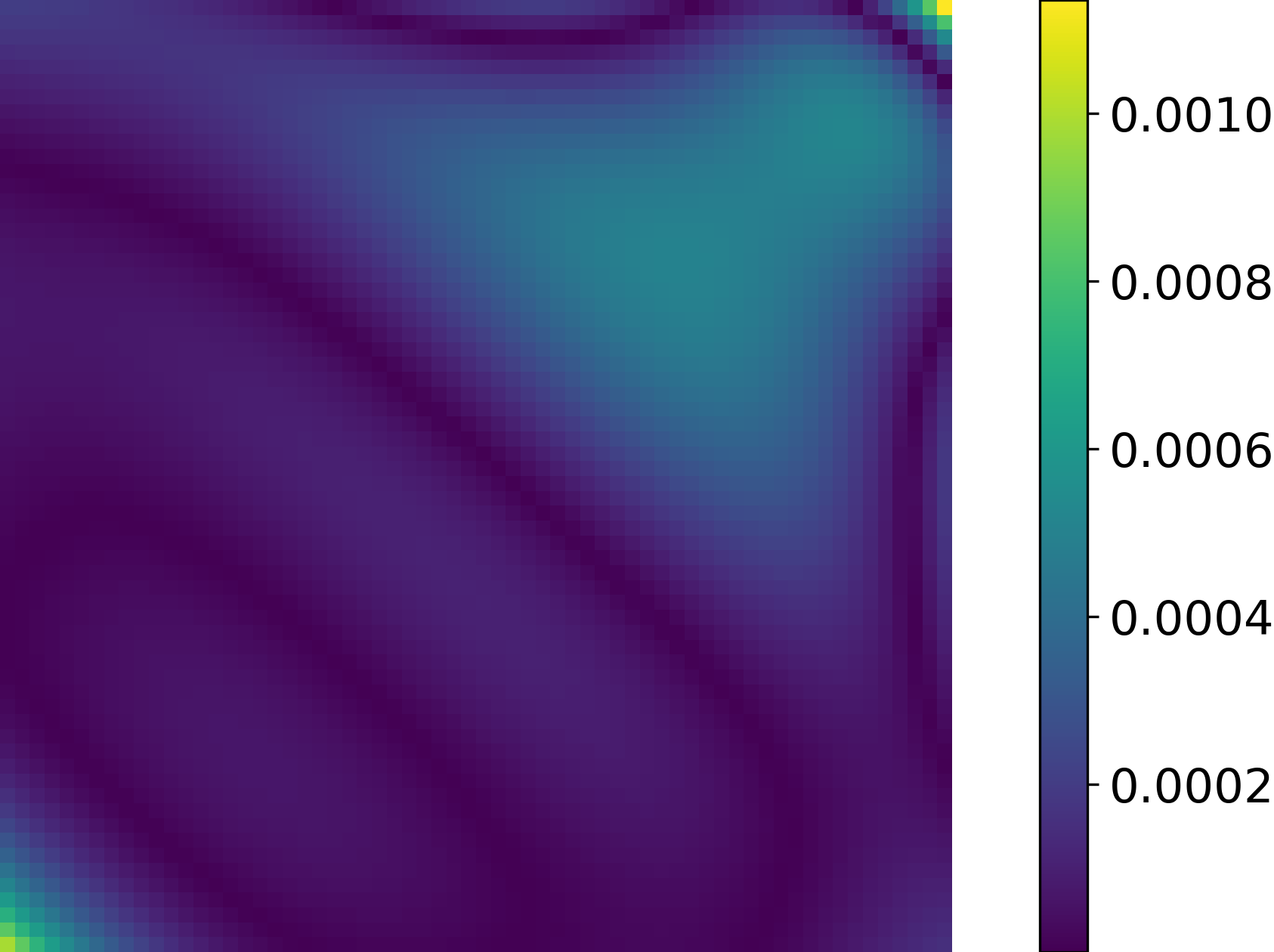}
        \vspace{-0.25cm} 
        \caption{absolute error $|k - \tilde{k}|$}
    \end{subfigure}

    \caption{Example 2-2. Visualizations of the exact coefficient $k$ and the reconstructed coefficient $\tilde{k}$ using 3D surface plots, along with the absolute error $|k - \tilde{k}|$ shown as a 2D heat map in the partially observed setting.}
    \label{fig:function_disconti_k}
\end{figure}
Figure~\ref{fig:function_disconti_k} presents the exact coefficient $k$, the reconstructed coefficient $\tilde{k}$, and the absolute error $|k - \tilde{k}|$, providing a visual comparison of the reconstruction performance under the partially observed setting.
The results confirm the robustness of the proposed method in accurately reconstructing the spatially varying coefficient, even with limited boundary measurements.

\begin{table}[!ht]
    \centering
    \begin{tabular}{|c|c|c|c|c|c|}
    \hline
    \textbf{Case} & \textbf{Observed} & \textbf{MSE} & \textbf{MAE} & \textbf{PSNR} & \textbf{SSIM} \\ \hline
    \multirow{2}{*}{Constant} 
        & Fully     & $3.10 \times 10^{-11}$ & $4.11 \times 10^{-6}$  & 105.09 & 0.999999999 \\ \cline{2-6}
        & Partially & $8.93 \times 10^{-10}$ & $2.11 \times 10^{-5}$  & 90.49  & 0.999999997 \\ \hline
    \multirow{2}{*}{Spatially varying} 
        & Fully     & $1.95 \times 10^{-8}$  & $1.11 \times 10^{-4}$  & 76.18  & 0.999999521 \\ \cline{2-6}
        & Partially & $4.10 \times 10^{-8}$  & $1.41 \times 10^{-4}$  & 72.95  & 0.999999575 \\ \hline
    \end{tabular}
    \caption{Example 1-1 - Example 2-2. Comparison of MSE, MAE, PSNR, and SSIM for the constant and spatially varying cases under fully and partially observed settings.}
    \label{tab:constant_function_metrics}
\end{table}
Table~\ref{tab:constant_function_metrics} presents the performance metrics—mean squared error (MSE), mean absolute error (MAE), peak signal-to-noise ratio (PSNR), and structural similarity index (SSIM)—for Examples~1-1, 1-2, 2-1, and 2-2.
It compares the reconstruction quality under fully and partially observed settings for both the constant and spatially varying coefficient cases.
PSNR measures the numerical similarity of the reconstructed values to the exact values~\cite{Huynh08}, while SSIM evaluates the structural consistency between the reconstructed and exact values~\cite{Wang04}.
Higher values of PSNR and SSIM indicate better reconstruction performance, with SSIM values closer to 1 reflecting greater structural similarity. In contrast, lower values of MSE and MAE correspond to smaller reconstruction errors.

In the constant case, the fully observed setting outperforms the partially observed setting in all metrics, including MSE, MAE, PSNR, and SSIM.
Nevertheless, the partially observed setting achieves a very low MSE of $8.93 \times 10^{-10}$ and a SSIM value close to 1, indicating that the overall reconstruction remains highly accurate.
In the spatially varying case, all metrics are slightly lower than in the constant case, reflecting the increased complexity of the coefficient.
The fully observed setting yields lower MSE and MAE values and achieves a higher PSNR, indicating a more accurate reconstruction.
Although the SSIM value is slightly higher in the partially observed setting, the difference is negligible.

\subsection{Example 3: Phantom Coefficient with a Zero Source Function}
\label{subsec:phantom}
In this subsection, we consider the phantom coefficient, which is used in many medical imaging applications such as EIT~\cite{Bar21, Cheney99, Patterson05}.

In this case, we only consider the fully observed setting, where $\omega = \partial\Omega$. The source function $f$ is set to zero, and the domain $\Omega$ is defined as a circular region:  
$\Omega = \{(x,y) \mid (x-0.5)^2 + (y-0.5)^2 \leq 0.5^2\}.$
The phantom coefficient $k$ within this domain is configured as shown in Figures~\ref{fig:phantom_1_no_k}(a) and~\ref{fig:phantom_2_no_k}(a).
We conduct numerical experiments on two different phantom coefficients $k$, where a slight smoothing effect is applied using a Gaussian kernel.
This reduces sharp transitions in the values and yields a more continuous distribution across the domain.

To enhance the accuracy in the case where the coefficient varies over the space, we include a regularization term in the loss function to improve the stability and robustness of the reconstructed coefficient~\cite{Bar21}, and use the following loss function:
\begin{equation}\label{eq:loss_with_gradient}
    \mathfrak{L}(\tilde{k},\mathfrak{D}, \tilde{\mathfrak{D}}) 
    = \sum_{i=1}^N \bigl\|h_{g_i} - \tilde{k}\,\partial_{\nu}\tilde{p}_{g_i}\bigr\|_{L^2(\partial\Omega)}^2 
    + \lambda \bigl\|\nabla \tilde{k}\bigr\|_{L^2(\Omega)}^2.
\end{equation}
Here, $\mathfrak{D} = \{(g_i, h_{g_i})_{i = 1, \dots, N} \mid g_i \in H^{1/2}(\partial\Omega), \, \operatorname{supp}g_i \subset \partial\Omega, h_{g_i} = k \partial_\nu p_{g_i}|_{\partial\Omega}\}$,  
$\tilde{\mathfrak{D}} = \{(g_i, \tilde{k}\partial_\nu \tilde{p}_{g_i}|_{\partial\Omega})_{i = 1, \dots, N} \mid g_i \in H^{1/2}(\partial\Omega), \, \operatorname{supp}g_i \subset \partial\Omega\}$.
The positive parameter $\lambda$ acts as a regularization weight, controlling the contribution of the regularization term $\|\nabla\tilde{k}\|_{L^2(\Omega)}^2$.

The experimental setup for the phantom case is as follows.
The Dirichlet boundary condition is imposed, given by $g(x,y) = \cos(a\pi x)\cos(b\pi y)\big|_{\partial\Omega},$
where $a$ and $b$ are uniformly sampled from the interval $[0,2]$.
The circular domain is discretized using a triangular mesh with 51,433 degrees of freedom for the linear Lagrange finite element space.
This experiment uses a Cauchy dataset containing 2,048 samples.
The learning rate is set to 0.001 for the first 50 epochs and reduced to 0.0001 for the remaining 50 epochs.

\subsubsection{Example 3-1-1: Phantom 1}

\begin{figure}[!ht]
    \centering
    \begin{subfigure}{0.3\textwidth}
        \centering
        \includegraphics[width=\textwidth]{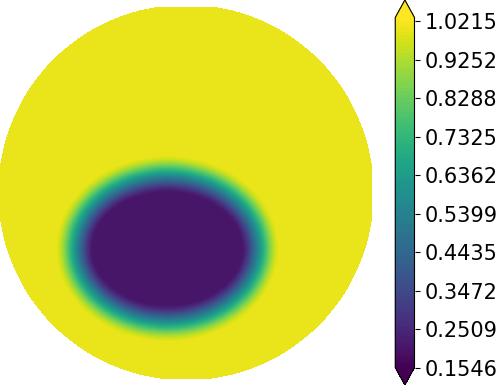}
        \caption{Phantom 1: exact $k$}
    \end{subfigure}
    \hfill
    \begin{subfigure}{0.3\textwidth}
        \centering
        \includegraphics[width=\textwidth]{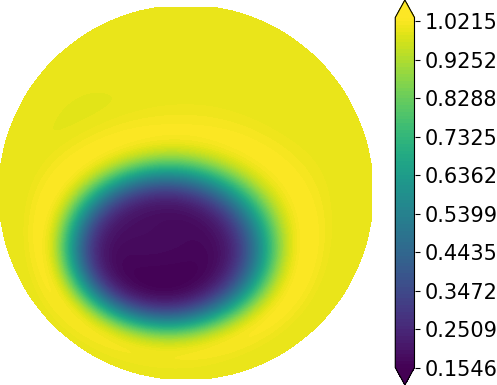}
        \caption{$\tilde{k}$ using~\eqref{eq:loss}}
    \end{subfigure}
    \hfill
    \begin{subfigure}{0.3\textwidth}
        \centering
        \includegraphics[width=\textwidth]{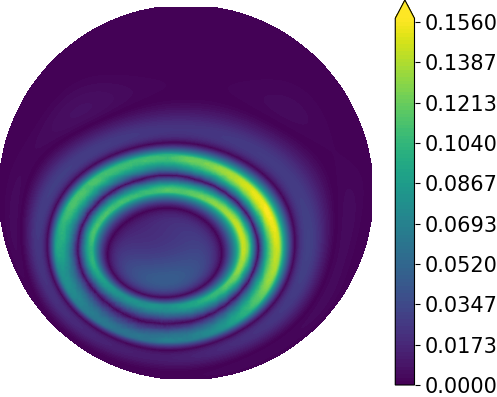}
        \caption{absolute error $|k - \tilde{k}|$}
    \end{subfigure}

    \hspace*{0.335\textwidth}
    \begin{subfigure}{0.3\textwidth}
        \includegraphics[width=\textwidth]{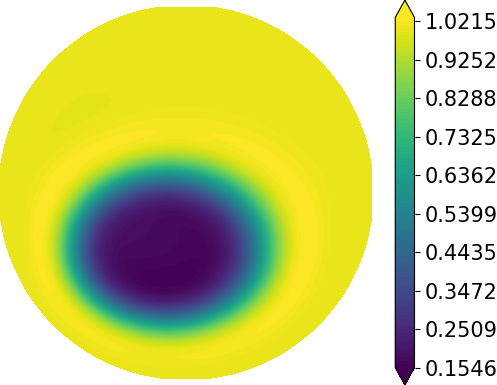}
        \caption{$\tilde{k}$ using~\eqref{eq:loss_with_gradient}}
    \end{subfigure}
    \hfill
    \begin{subfigure}{0.3\textwidth}
        \centering
        \includegraphics[width=\textwidth]{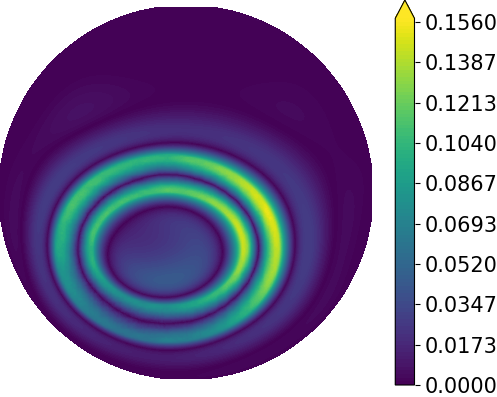}
        \caption{absolute error $|k - \tilde{k}|$}
    \end{subfigure}

    \caption{Example 3-1-1. Phantom 1. Visualizations of the exact coefficient $k$, reconstructed coefficient $\tilde{k}$, and the corresponding absolute errors under two different loss functions.
     (c) shows the absolute error $|k - \tilde{k}|$ between (a) and (b), while (e) shows the absolute error between (a) and (d).}
    \label{fig:phantom_1_no_k}
\end{figure}

In this case, Phantom 1 (Figure~\ref{fig:phantom_1_no_k}(a)) is used as the coefficient $k$ due to its relatively simple structure compared to the other phantom.
The regularization weight in the loss function~\eqref{eq:loss_with_gradient} is set to $\lambda = 2 \times 10^{-8}$.

\begin{figure}[!ht]
    \centering
    \begin{subfigure}{0.3\textwidth}
        \centering
        \includegraphics[width=\textwidth]{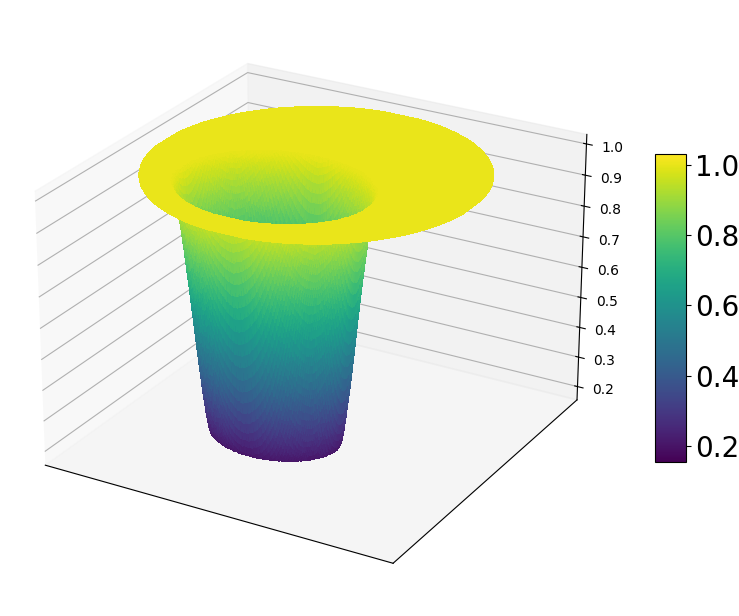}
        \caption{exact $k$}
    \end{subfigure}
    \hfill
    \begin{subfigure}{0.3\textwidth}
        \centering
        \includegraphics[width=\textwidth]{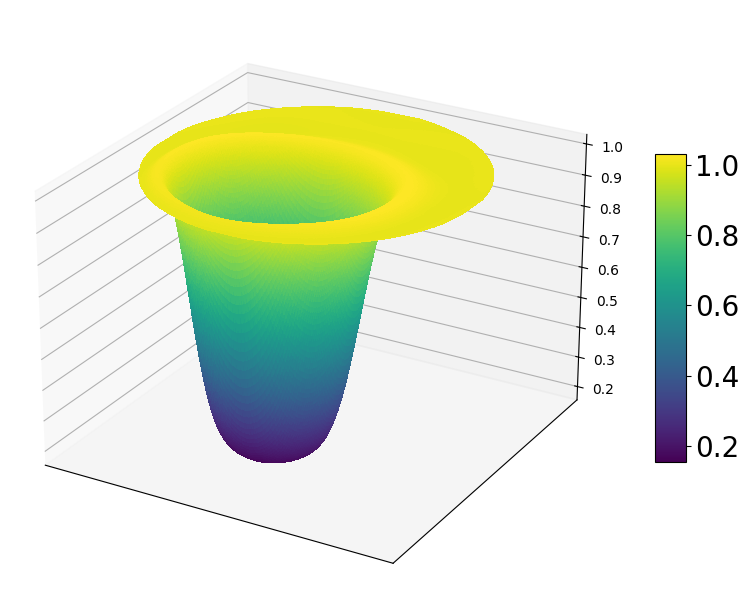}
        \caption{$\tilde{k}$ using~\eqref{eq:loss}}
    \end{subfigure}
    \hfill
    \begin{subfigure}{0.3\textwidth}
        \centering
        \includegraphics[width=\textwidth]{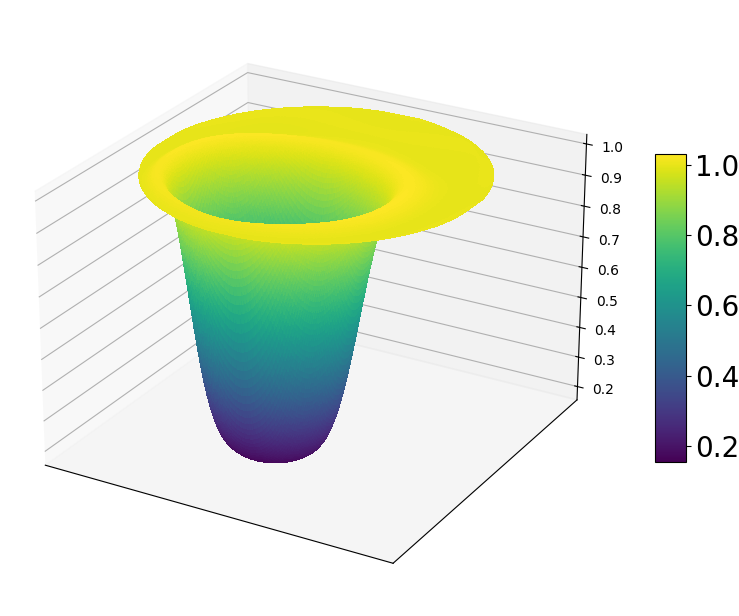}
        \caption{$\tilde{k}$ using~\eqref{eq:loss_with_gradient}}
    \end{subfigure}

    \caption{Example 3-1-1. Three-dimensional visualizations of the exact coefficient $k$ and the reconstructed coefficient $\tilde{k}$ obtained using the loss functions~\eqref{eq:loss} and~\eqref{eq:loss_with_gradient}.}
    \label{fig:phantom_1_no_k_3d}
\end{figure}

Figure~\ref{fig:phantom_1_no_k} compares the results of two models: one trained with the regularization term included in the loss function~\eqref{eq:loss_with_gradient}, and the other trained without it using the loss function~\eqref{eq:loss}.
The results indicate that the inclusion of the regularization term does not lead to a significant difference in performance.
Additionally, Figure~\ref{fig:phantom_1_no_k_3d} provides three-dimensional visualizations of the phantom coefficient.

\begin{figure}[ht!]
    \centering
    \includegraphics[width=0.6\textwidth]{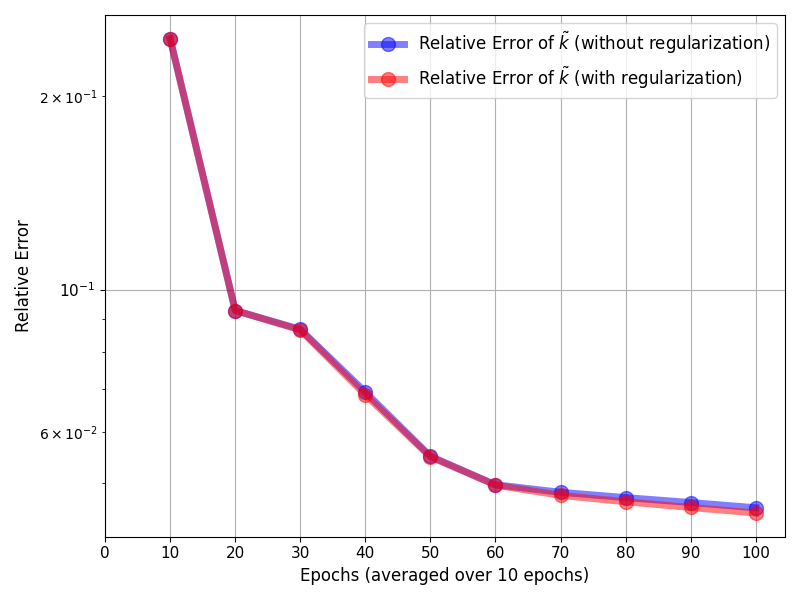}
    \caption{Example 3-1-1. Relative error of $\tilde{k}$ with and without regularization over epochs (log scale).  Each point represents the averaged value over a 10-epoch interval.}
    \label{fig:phantom_1_no_combined}
\end{figure}
Figure~\ref{fig:phantom_1_no_combined} shows the relative error of $\tilde{k}$ on a logarithmic scale.
The results indicate that the model trained with the regularization slightly outperforms the model trained without regularization, although the difference remains negligible.
Both models exhibit a steady decrease in relative error and eventually converge to similar values.
This suggests that, in this case, the regularization has little impact on the reconstruction accuracy of the coefficient $k$.

\subsubsection{Example 3-1-2: Phantom 1 with Post-Processing}
In this subsection, we extend the experiment conducted on Phantom 1 in Example 3-1-1 by incorporating a post-processing step.
In this step, we apply the hyperbolic tangent function to constrain the reconstructed coefficient $\tilde{k}$ within the range $[0.2, 1.0]$, based on the assumed prior knowledge of the range of $k$.
The regularization weight in the loss function~\eqref{eq:loss_with_gradient} is set to $\lambda = 2 \times 10^{-7}$.
All other experimental configurations remain the same as in Example 3-1-1.

\begin{figure}[!ht]
    \centering
    \begin{subfigure}{0.3\textwidth}
        \centering
        \includegraphics[width=\textwidth]{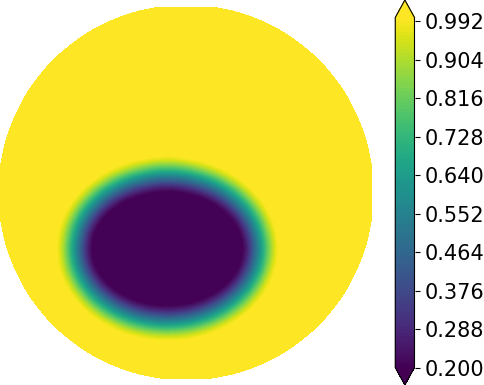}
        \caption{Phantom 1: exact $k$}
    \end{subfigure}
    \hfill
    \begin{subfigure}{0.3\textwidth}
        \centering
        \includegraphics[width=\textwidth]{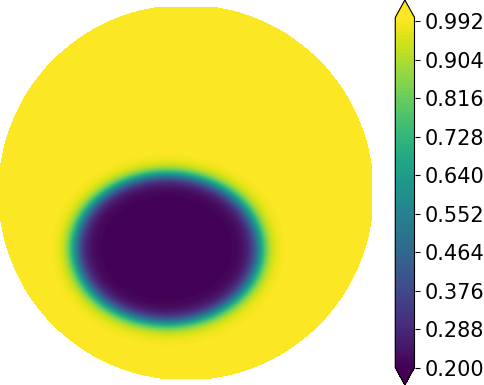}
        \caption{$\tilde{k}$ using~\eqref{eq:loss}}
    \end{subfigure}
    \hfill
    \begin{subfigure}{0.3\textwidth}
        \centering
        \includegraphics[width=\textwidth]{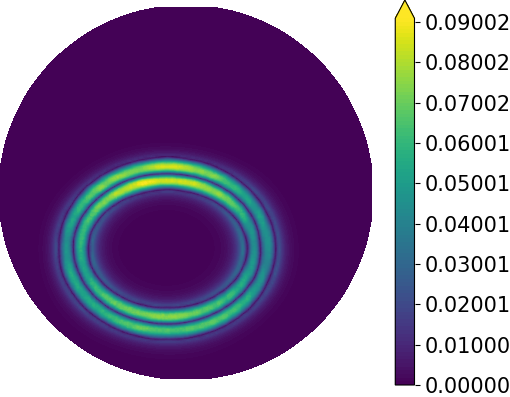}
        \caption{absolute error $|k - \tilde{k}|$}
    \end{subfigure}

    \hspace*{0.335\textwidth}
    \begin{subfigure}{0.3\textwidth}
        \includegraphics[width=\textwidth]{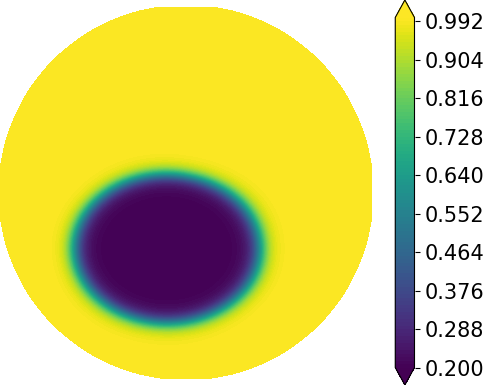}
        \caption{$\tilde{k}$ using~\eqref{eq:loss_with_gradient}}
    \end{subfigure}
    \hfill
    \begin{subfigure}{0.3\textwidth}
        \centering
        \includegraphics[width=\textwidth]{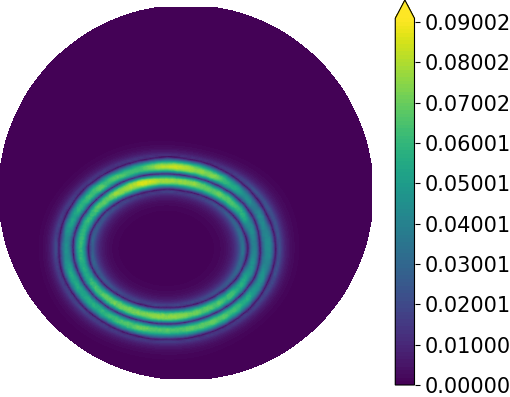}
        \caption{absolute error $|k - \tilde{k}|$}
    \end{subfigure}

    \caption{Example 3-1-2. Visualizations of the exact coefficient $k$, reconstructed coefficient $\tilde{k}$, and the corresponding absolute errors under two different loss functions with the proposed post-processing.
     (c) shows the absolute error $|k - \tilde{k}|$ between (a) and (b), while (e) shows the absolute error between (a) and (d).}
    \label{fig:phantom_1_k}
\end{figure}
Figure~\ref{fig:phantom_1_k} compares the results of models trained with and without the regularization term, using the loss functions~\eqref{eq:loss_with_gradient} and~\eqref{eq:loss}, respectively.
The comparison indicates that the inclusion of the regularization term does not significantly affect the final reconstruction, as both models produce similar outcomes.
In contrast, a comparison between the results without post-processing (Figure~\ref{fig:phantom_1_no_k}) and with post-processing (Figure~\ref{fig:phantom_1_k})  shows a significant reduction in the region with noticeable error after post-processing.
This suggests that the post-processing step plays an important role in improving the spatial accuracy of the reconstructed coefficient.

\begin{figure}[!ht]
    \centering
    \includegraphics[width=0.6\textwidth]{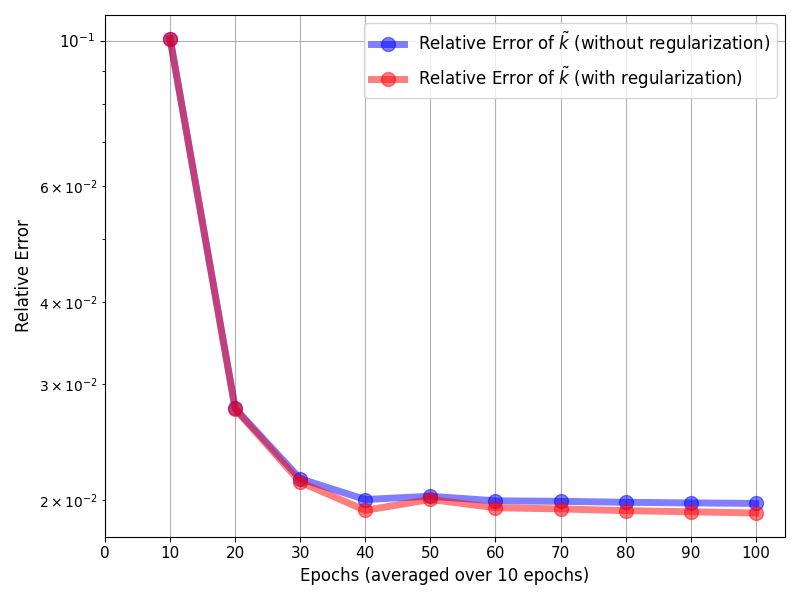}
    \caption{Example 3-1-2. Relative error of $\tilde{k}$ with and without regularization terms, using the post-processing, over epochs (log scale). Each point represents the average over a 10-epoch interval.}
    \label{fig:phantom_1_k_error}
\end{figure}
Figure~\ref{fig:phantom_1_k_error} presents the relative error in the coefficient $k$ for models trained with and without the regularization term, after applying the post-processing step.
As in the case without post-processing, the model trained with the regularization term shows slightly better performance, although the difference remains negligible.

\subsubsection{Example 3-2-1: Phantom 2}

\begin{figure}[!ht]
    \centering
    \begin{subfigure}{0.3\textwidth}
        \centering
        \includegraphics[width=\textwidth]{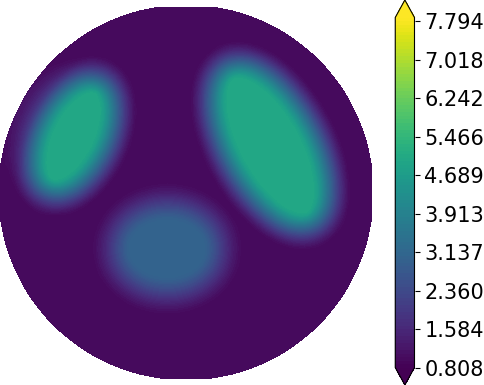}
        \caption{Phantom 2: exact $k$}
    \end{subfigure}
    \hfill
    \begin{subfigure}{0.3\textwidth}
        \centering
        \includegraphics[width=\textwidth]{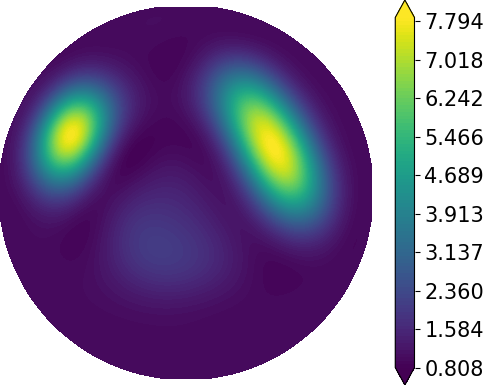}
        \caption{$\tilde{k}$ using~\eqref{eq:loss}}
    \end{subfigure}
    \hfill
    \begin{subfigure}{0.3\textwidth}
        \centering
        \includegraphics[width=\textwidth]{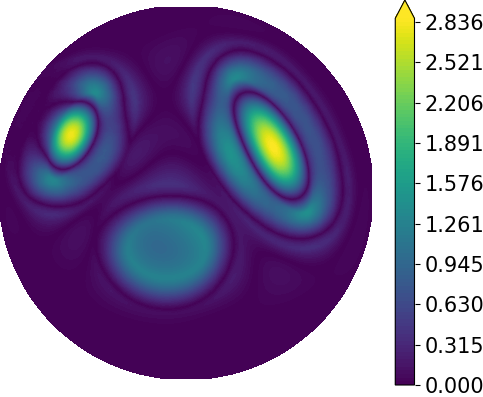}
        \caption{absolute error $|k - \tilde{k}|$}
    \end{subfigure}

    \hspace*{0.335\textwidth}
    \begin{subfigure}{0.3\textwidth}
        \centering
        \includegraphics[width=\textwidth]{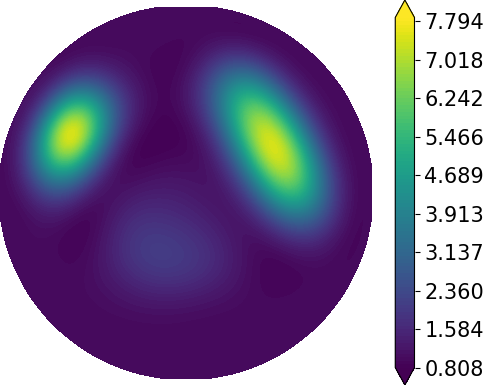}
        \caption{$\tilde{k}$ using~\eqref{eq:loss_with_gradient}}
    \end{subfigure}
    \hfill
    \begin{subfigure}{0.3\textwidth}
        \centering
        \includegraphics[width=\textwidth]{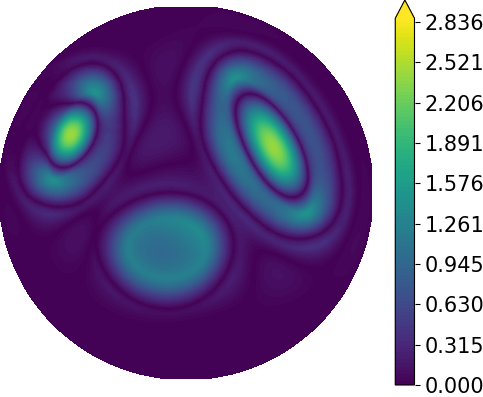}
        \caption{absolute error $|k - \tilde{k}|$}
    \end{subfigure}

    \caption{Example 3-2-1. Phantom 2. Visualizations of the exact coefficient $k$, reconstructed coefficient $\tilde{k}$, and the corresponding absolute errors under two different loss functions.
     (c) shows the absolute error $|k - \tilde{k}|$ between (a) and (b), while (e) shows the absolute error between (a) and (d).}
    \label{fig:phantom_2_no_k}
\end{figure}

In this case, a more complex phantom, Phantom 2 (Figure~\ref{fig:phantom_2_no_k}(a)), is used to more clearly assess the effect of the regularization term.
The regularization weight in the loss function~\eqref{eq:loss_with_gradient} is set to $\lambda = 2 \times 10^{-6}$.

Figure~\ref{fig:phantom_2_no_k} presents visualizations of the exact coefficient $k$, the reconstructed coefficient $\tilde{k}$ obtained using the loss functions~\eqref{eq:loss} and~\eqref{eq:loss_with_gradient}, and their corresponding absolute errors.
\begin{figure}[!ht]
    \centering
    \begin{subfigure}{0.3\textwidth}
        \centering
        \includegraphics[width=\textwidth]{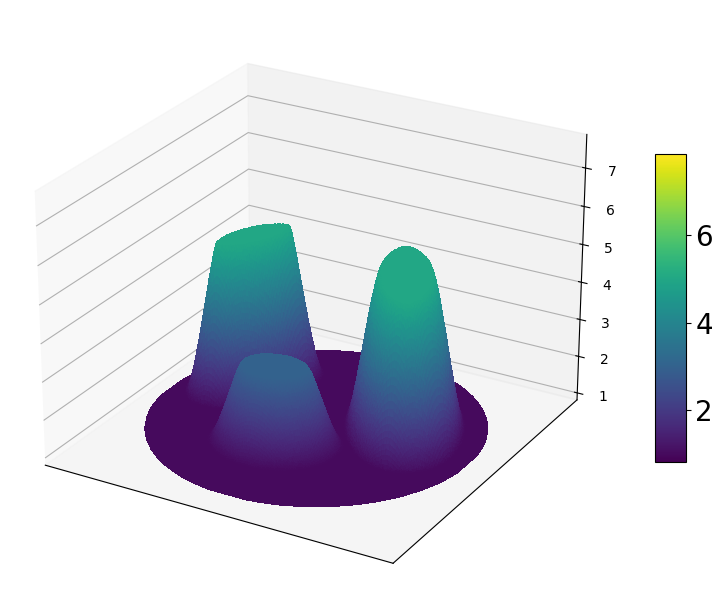}
        \caption{exact $k$}
    \end{subfigure}
    \hfill
    \begin{subfigure}{0.3\textwidth}
        \centering
        \includegraphics[width=\textwidth]{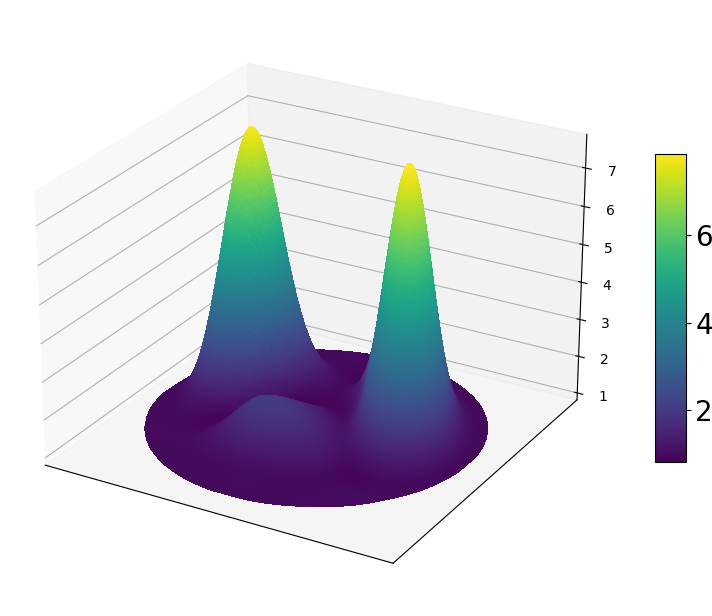}
        \caption{$\tilde{k}$ using~\eqref{eq:loss}}
    \end{subfigure}
    \hfill
    \begin{subfigure}{0.3\textwidth}
        \centering
        \includegraphics[width=\textwidth]{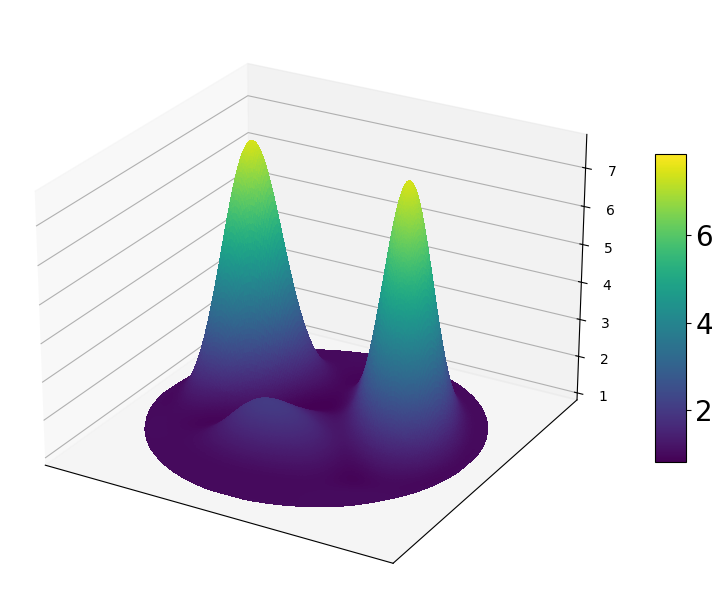}
        \caption{$\tilde{k}$ using~\eqref{eq:loss_with_gradient}}
    \end{subfigure}

   \caption{Example 3-2-1. Three-dimensional visualizations of the exact coefficient $k$ and the reconstructed coefficient $\tilde{k}$ obtained using the loss functions~\eqref{eq:loss} and~\eqref{eq:loss_with_gradient}.}
    \label{fig:phantom_2_no_k_3d}
\end{figure}
Furthermore, Figure~\ref{fig:phantom_2_no_k_3d} provides three-dimensional visualizations of the same results.
The overall structure of $k$ is well captured, but the reconstruction tends to overshoot the peaks, producing sharper and taller spikes than the exact coefficient.

\begin{figure}[ht!]
    \centering
    \includegraphics[width=0.6\textwidth]{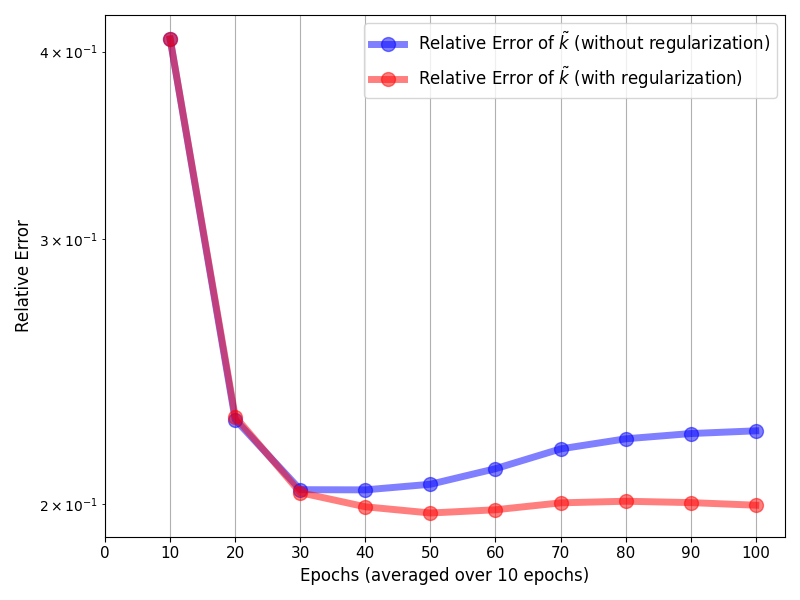}
    \caption{Example 3-2-1. Relative error of $\tilde{k}$ with and without regularization terms over epochs (log scale).  Each point represents the averaged value over a 10-epoch interval.}
    \label{fig:phantom_2_no_combined}
\end{figure}
Figure~\ref{fig:phantom_2_no_combined} shows the relative error of $\tilde{k}$ over epochs on a logarithmic scale.
While both models show similar convergence in the initial training epochs, the model trained without regularization shows a slight increase in error in the later stages.
In contrast, the model with regularization maintains stable convergence throughout, suggesting that the regularization term contributes to improved training stability.

\subsubsection{Example 3-2-2: Phantom 2 with Post-Processing}
In this subsection, we extend the experiment conducted on Phantom 2 in Example 3-2-1 by incorporating a post-processing step.
In this step, we apply the hyperbolic tangent function to constrain the reconstructed coefficient $\tilde{k}$ within the range $[1.0, 5.0]$, based on the assumed prior knowledge of the range of $k$.
The regularization weight in the loss function~\eqref{eq:loss_with_gradient} is set to $\lambda = 2 \times 10^{-6}$.
All other experimental configurations remain the same as in Example 3-2-1.

\begin{figure}[!ht]
    \centering
    \begin{subfigure}{0.3\textwidth}
        \centering
        \includegraphics[width=\textwidth]{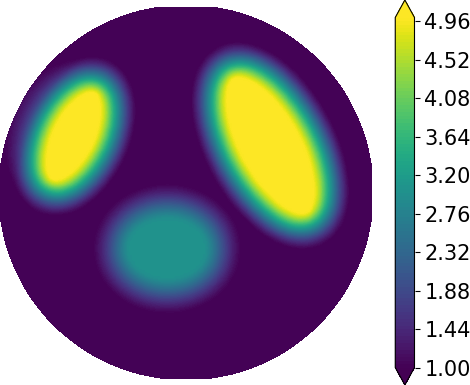}
        \caption{Phantom 2: exact $k$}
    \end{subfigure}
    \hfill
    \begin{subfigure}{0.3\textwidth}
        \centering
        \includegraphics[width=\textwidth]{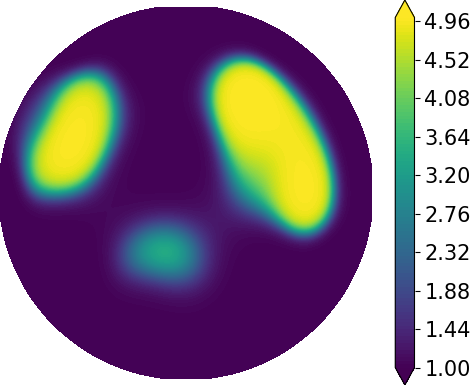}
        \caption{$\tilde{k}$ using~\eqref{eq:loss}}
    \end{subfigure}
    \hfill
    \begin{subfigure}{0.3\textwidth}
        \centering
        \includegraphics[width=\textwidth]{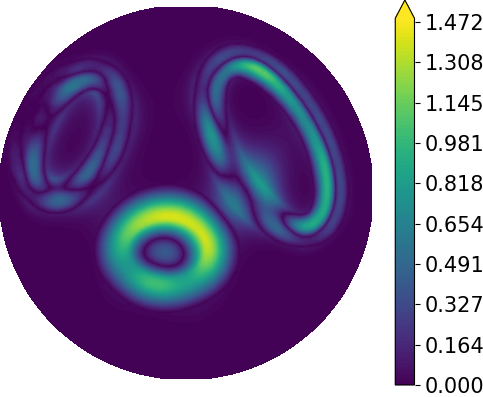}
        \caption{absolute error $|k - \tilde{k}|$}
    \end{subfigure}

    \hspace*{0.335\textwidth}
    \begin{subfigure}{0.3\textwidth}
        \includegraphics[width=\textwidth]{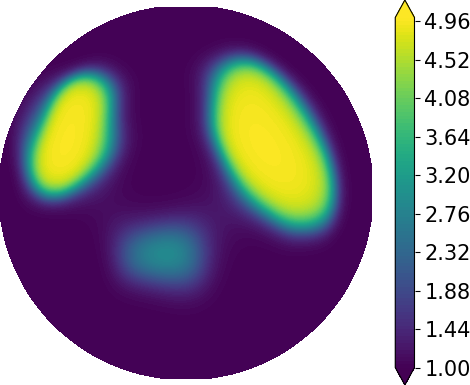}
        \caption{$\tilde{k}$ using~\eqref{eq:loss_with_gradient}}
    \end{subfigure}
    \hfill
    \begin{subfigure}{0.3\textwidth}
        \centering
        \includegraphics[width=\textwidth]{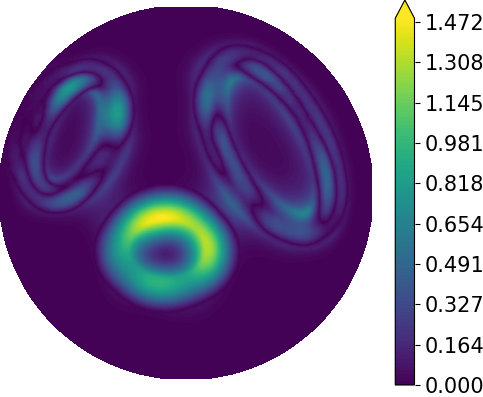}
        \caption{absolute error $|k - \tilde{k}|$}
    \end{subfigure}

    \caption{Example 3-2-2.Visualizations of the exact coefficient $k$, reconstructed coefficient $\tilde{k}$, and the corresponding absolute errors under two different loss functions with post-processing.
     (c) shows the absolute error $|k - \tilde{k}|$ between (a) and (b), while (e) shows the absolute error between (a) and (d).}
    \label{fig:phantom_2_k}
\end{figure}
Figure~\ref{fig:phantom_2_k} compares the reconstruction results of two models: one trained using the loss function~\eqref{eq:loss}, and the other using the loss function~\eqref{eq:loss_with_gradient}, which includes a regularization term.
Figures~\ref{fig:phantom_2_k}(c) and~\ref{fig:phantom_2_k}(e) show the corresponding absolute error distributions.
In Figure~\ref{fig:phantom_2_k}(c), large errors are observed near regions with sharp transitions in $k$, while Figure~\ref{fig:phantom_2_k}(e) shows smaller errors in these regions.
This indicates that the model trained with the regularization term is more effective at capturing rapid variations and interfaces in the coefficient.
These results suggest that the regularization term plays a critical role in preserving the structural consistency of $k$, particularly in complex phantom coefficients.
\begin{figure}[!ht]
    \centering
    \begin{subfigure}{0.3\textwidth}
        \centering
        \includegraphics[width=\textwidth]{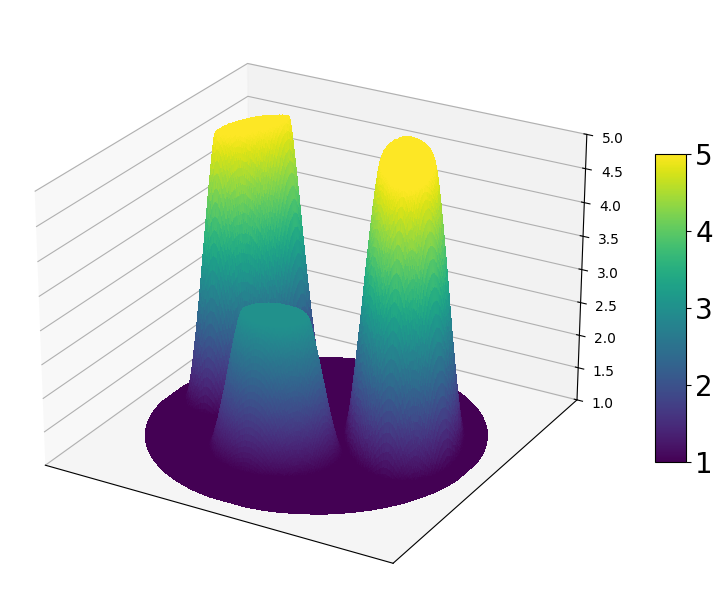}
        \caption{exact $k$}
    \end{subfigure}
    \hfill
    \begin{subfigure}{0.3\textwidth}
        \centering
        \includegraphics[width=\textwidth]{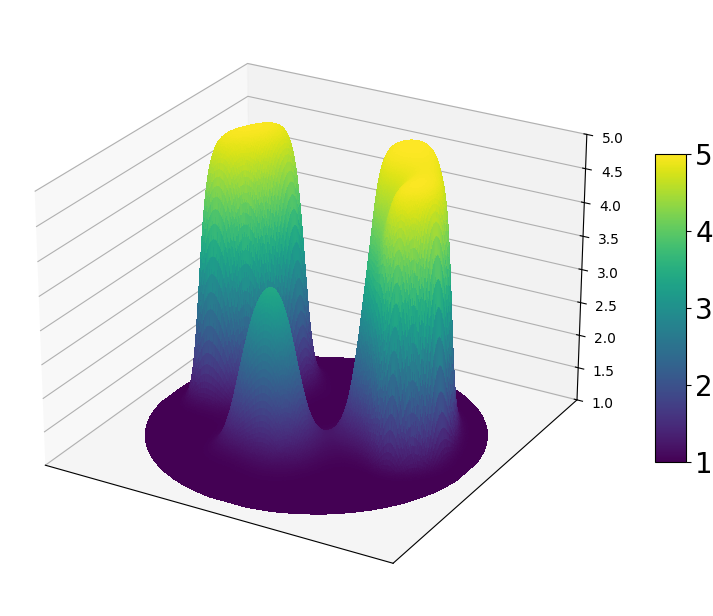}
        \caption{$\tilde{k}$ using~\eqref{eq:loss}}
    \end{subfigure}
    \hfill
    \begin{subfigure}{0.3\textwidth}
        \centering
        \includegraphics[width=\textwidth]{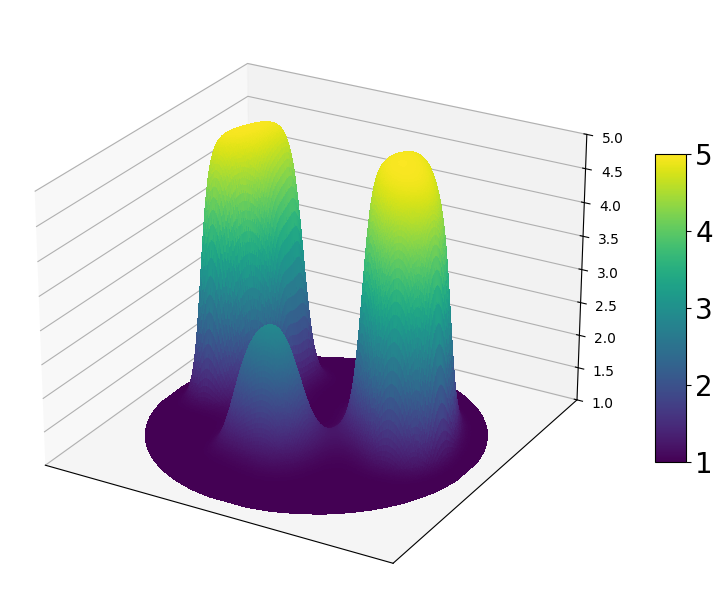}
        \caption{$\tilde{k}$ using~\eqref{eq:loss_with_gradient}}
    \end{subfigure}

    \caption{Example 3-2-2. Three-dimensional visualizations of the exact coefficient $k$ and the reconstructed coefficient $\tilde{k}$ obtained using the loss functions~\eqref{eq:loss} and~\eqref{eq:loss_with_gradient} with post-processing.}
    \label{fig:phantom_2_k_3d}
    
\end{figure}
To complement this, Figure~\ref{fig:phantom_2_k_3d} presents three-dimensional visualizations of the same results.
These visualizations provide a clearer perspective on the shape and height of the reconstructed peaks.
While both models capture the overall structure of $k$, the model with the regularization term—shown in Figure~\ref{fig:phantom_2_k_3d}(c)—achieves more accurate peak heights and smoother transitions.
In contrast, the model without the regularization term, shown in Figure~\ref{fig:phantom_2_k_3d}(b), shows slight deformation near regions with sharp transitions in $k$.

\begin{figure}[!ht]
    \centering
    \includegraphics[width=0.6\textwidth]{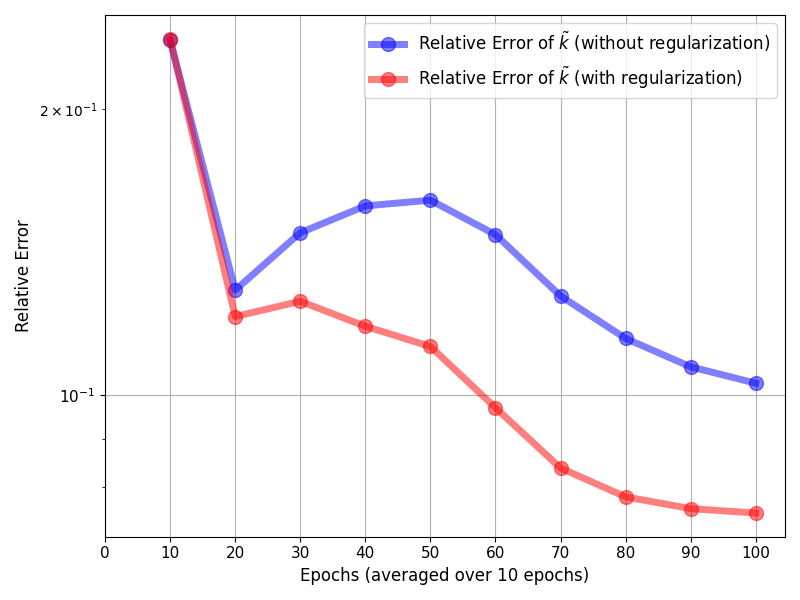}
    \caption{Example 3-2-2. Relative error of $\tilde{k}$ with and without regularization terms, using the post-processing, over epochs (log scale).  Each point represents the averaged value over a 10-epoch interval.}
    \label{fig:phantom_2_k_error}
\end{figure}
Figure~\ref{fig:phantom_2_k_error} shows the relative error of $\tilde{k}$ over epochs, plotted on a logarithmic scale.
The model with regularization (red curve) shows a consistently decreasing error, especially after 30 epochs, and ultimately reaches a significantly lower final error than the model without regularization (blue curve).
This behavior suggests that the regularization term not only improves reconstruction accuracy but also enhances training stability by promoting reliable convergence in complex scenarios.

\begin{table}[!ht]
    \centering
    \begin{tabular}{|c|c|c|c|c|c|c|}
    \hline
    \textbf{Case} & \textbf{Post-Processing} & \textbf{Regularization} & \textbf{MSE} & \textbf{MAE} & \textbf{PSNR} & \textbf{SSIM} \\ \hline
    \multirow{4}{*}{Phantom 1} 
        & \multirow{2}{*}{No}  & No  & 0.00162 & 0.02349 & 28.17 & 0.9575 \\ \cline{3-7}
                                &                          & Yes & 0.00155 & 0.02291 & 28.33 & 0.9583 \\ \cline{2-7}
        & \multirow{2}{*}{Yes} & No  & 0.00031 & 0.00721 & 35.03 & 0.9868 \\ \cline{3-7}
                                &                          & Yes & 0.00030 & 0.00708 & 35.30 & 0.9872 \\ \hline
    \multirow{4}{*}{Phantom 2} 
        & \multirow{2}{*}{No}  & No  & 0.26602 & 0.25938 & 23.66 & 0.9336 \\ \cline{3-7}
                                &                          & Yes & 0.21033 & 0.23995 & 24.20 & 0.9355 \\ \cline{2-7}
        & \multirow{2}{*}{Yes} & No  & 0.05494 & 0.11752 & 26.58 & 0.9488 \\ \cline{3-7}
                                &                          & Yes & 0.02955 & 0.09589 & 29.27 & 0.9653 \\ \hline
    \end{tabular}
    \caption{Example 3.1-1 - 3.2-2. Comparison of metrics—MSE, MAE, PSNR, and SSIM—for Phantom 1 and Phantom 2 under different settings of post-processing and regularization.}
    \label{tab:phantom_metrics}
\end{table}
Table~\ref{tab:phantom_metrics} summarizes the performance metrics—MSE, MAE, PSNR, and SSIM—for Phantom 1 and Phantom 2 under different combinations of post-processing and regularization.
For both phantoms, post-processing leads to significant improvements in all metrics.
These results suggest that while both regularization and post-processing improve reconstruction accuracy, their combined effect becomes particularly beneficial in more complex scenarios, such as Phantom~2.

\subsection{Example 4: Discontinuous Coefficient}

\begin{figure}[!ht]
    \centering
    \begin{subfigure}{0.3\textwidth}
        \centering
        \includegraphics[width=\textwidth]{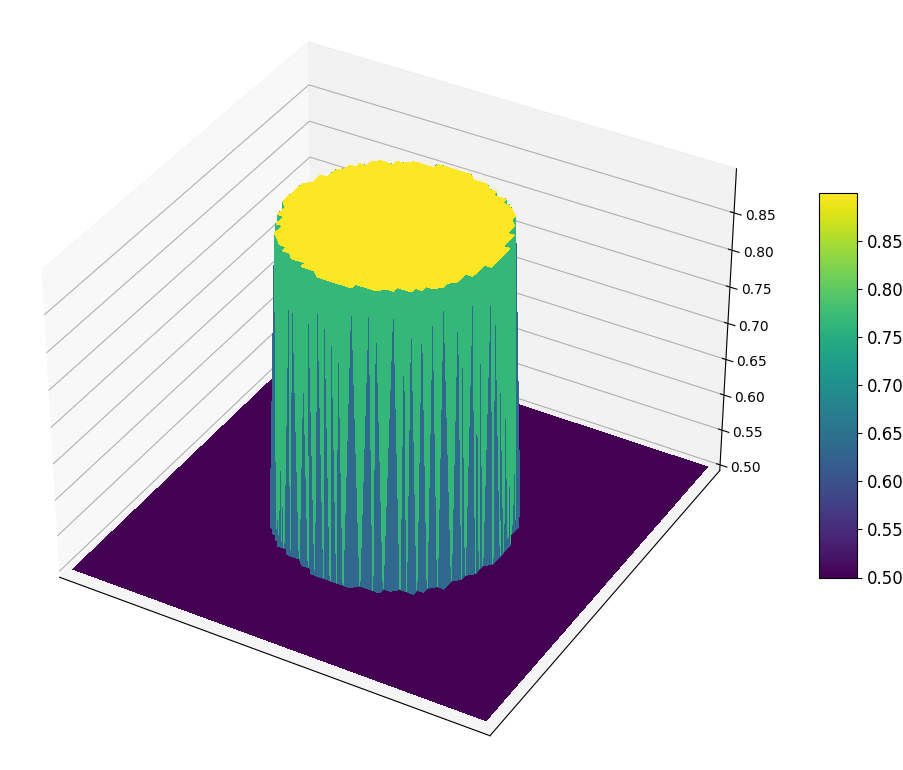}
        \caption{exact $k$}
    \end{subfigure}
    \hfill
    \begin{subfigure}{0.3\textwidth}
        \centering
        \includegraphics[width=\textwidth]{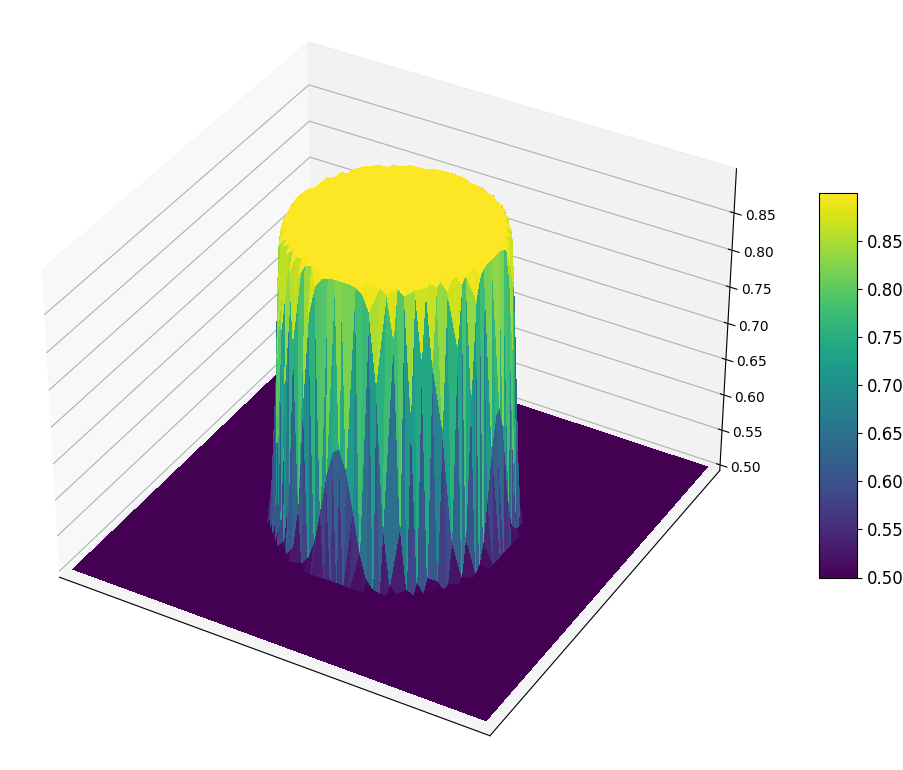}
        \caption{$\tilde{k}$ using~\eqref{eq:loss} with $\omega = \partial \Omega$}
    \end{subfigure}
    \hfill
    \begin{subfigure}{0.3\textwidth}
        \centering
        \includegraphics[width=\textwidth]{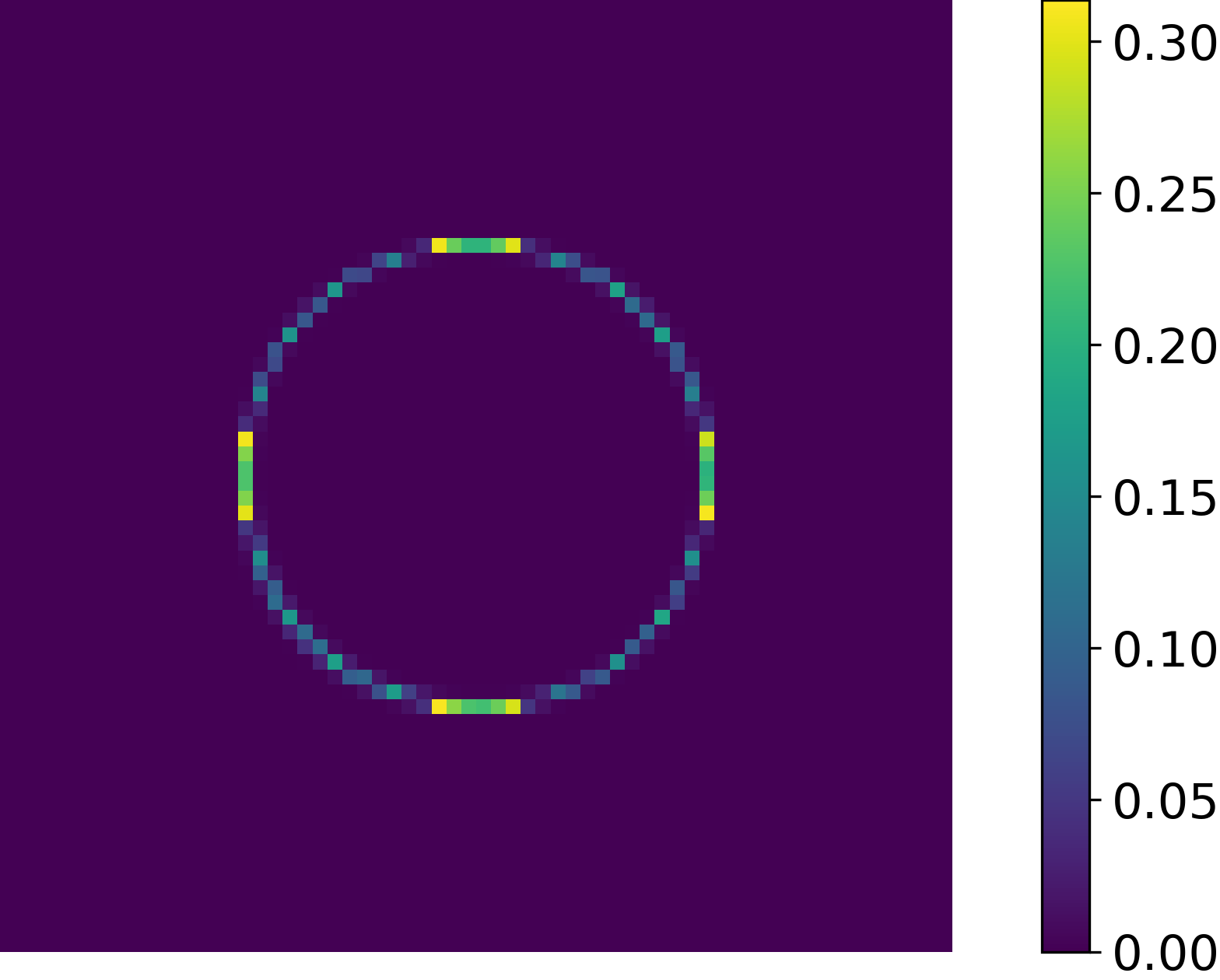} 
        \caption{absolute error $|k - \tilde{k}|$}
    \end{subfigure}

    \hspace*{0.335\textwidth}
    \begin{subfigure}{0.3\textwidth}
        \centering
        \includegraphics[width=\textwidth]{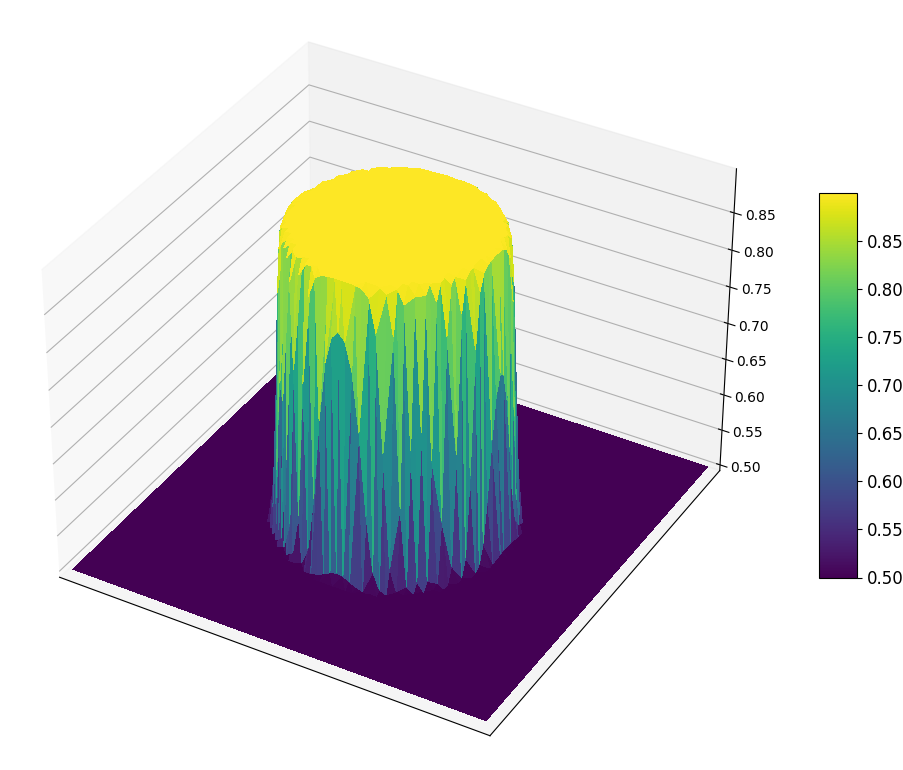}
        \caption{$\tilde{k}$ using~\eqref{eq:loss} with $\omega = \Gamma$}
    \end{subfigure}
    \hfill
    \begin{subfigure}{0.3\textwidth}
        \centering
        \includegraphics[width=\textwidth]{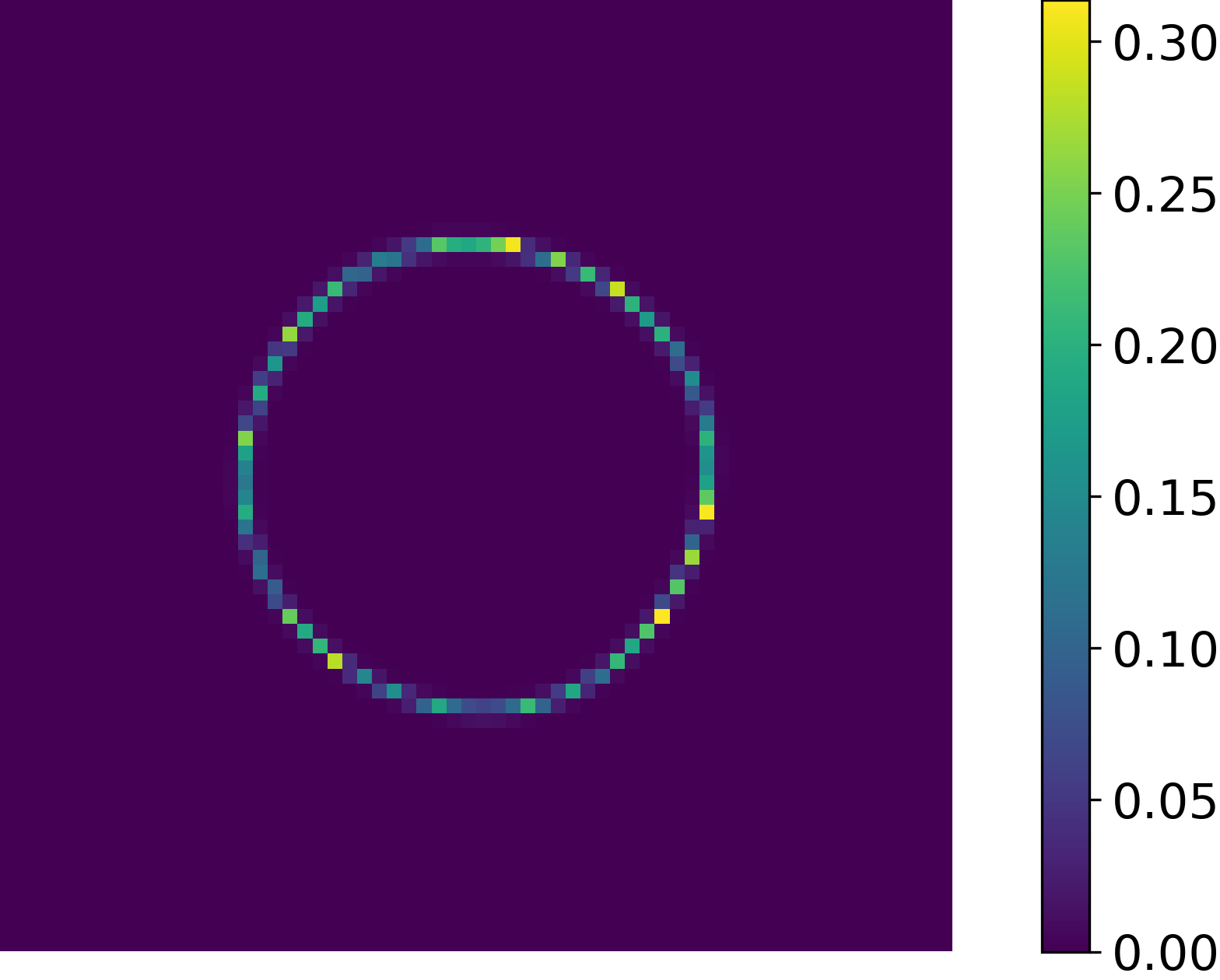}
        \caption{absolute error $|k - \tilde{k}|$}
    \end{subfigure}
    \caption{Example 4.  Visualizations of the exact coefficient $k$, the reconstructed coefficient $\tilde{k}$ in three dimensions, and the corresponding absolute errors as a heatmap, for two different observation settings.
        (c) shows the absolute error $|k - \tilde{k}|$ between (a) and (b), while (e) shows the error between (a) and (d).}    
    \label{fig:discontinuous_k}
\end{figure}

In this final example, we consider a more challenging and interesting case than the previous one, where the coefficient $k$ is discontinuous.
The coefficient is defined as
\begin{equation*}
k(x,y) = \left\{
\begin{array}{rl}
0.9 & \text{ if } (x-0.5)^2 + (y-0.5)^2 < (0.25)^2,\\
0.5 & \text{ otherwise},
\end{array}\right.
\end{equation*}
as illustrated in Figure~\ref{fig:discontinuous_k}(a).

We conduct experiments for both the fully and partially observed settings.
All other configurations, including the computational domain, degrees of freedom, Dirichlet boundary condition, number of data samples, and learning rate, are the same as in Example 1-2.
The source function is defined as  $f(x,y) = -\nabla \cdot (k(x,y) \nabla (\cos(a\pi x)\cos(b\pi y)))$.
For the partially observed setting, this source function is imposed on $\Gamma$, while $f$ is set to zero on $\Gamma_0$.
We also include the post-processing step by applying the hyperbolic tangent function to constrain the reconstructed coefficient $\tilde{k}$ within the range $[0.4, 1.0]$, based on the assumed prior knowledge of the range of $k$.
In this example, we train the model using the loss function~\eqref{eq:loss}, considering both $\omega = \partial\Omega$ and $\omega = \Gamma$.

Figure~\ref{fig:discontinuous_k} shows the reconstruction results and the corresponding absolute errors for both the fully observed and partially observed settings.
The absolute error maps in Figures~\ref{fig:discontinuous_k}(c) and~\ref{fig:discontinuous_k}(e) show that the reconstruction errors are primarily concentrated near regions where $k$ is discontinuous.
This indicates that the model has difficulty capturing sharp transitions, which is expected due to the inherent challenge of learning discontinuous coefficients.

\begin{figure}[!ht]
    \centering
    \begin{subfigure}{0.7\textwidth}
        \centering
        \includegraphics[width=\textwidth]{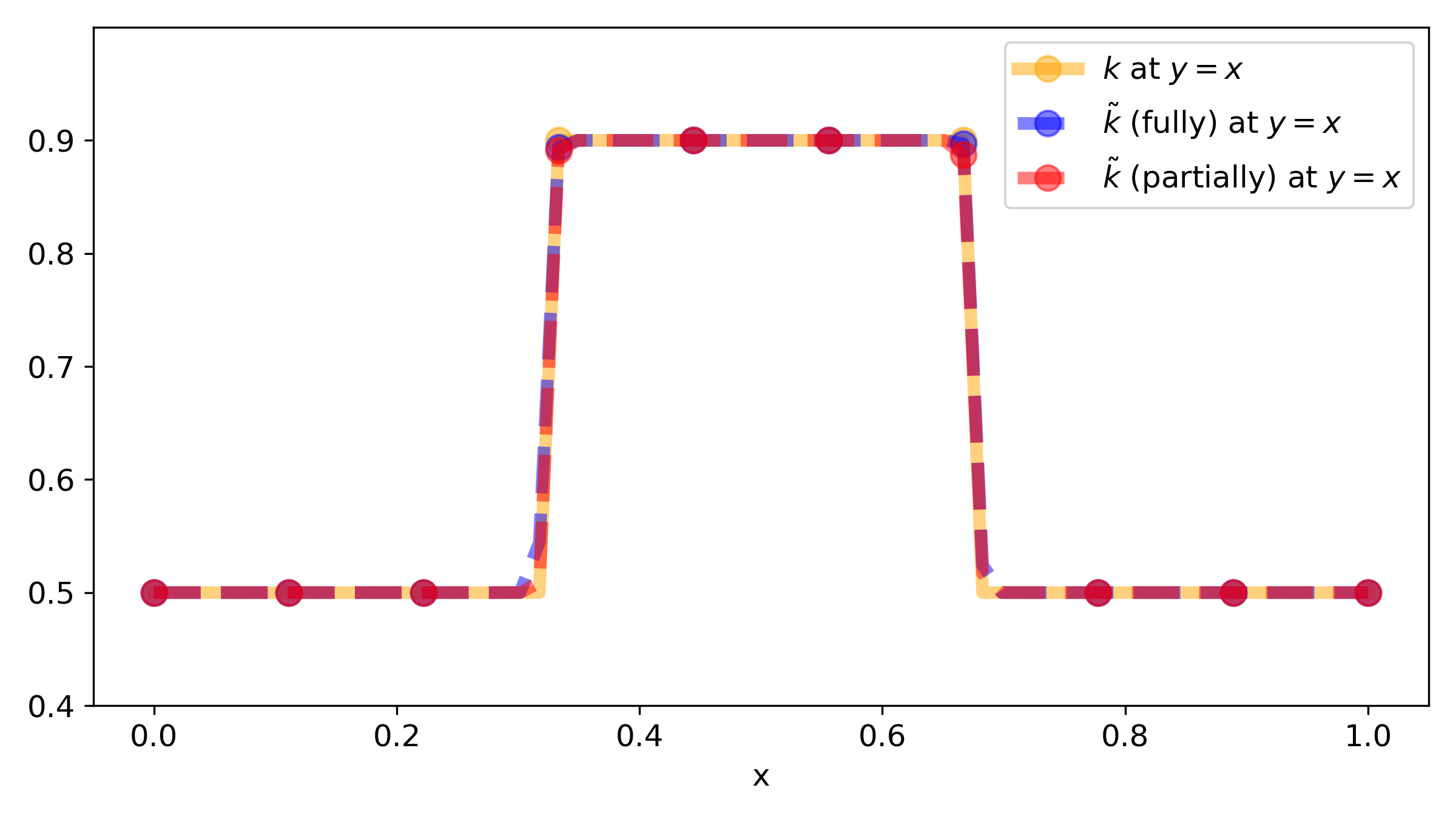}
        \caption{$y = x$}
    \end{subfigure}
    \caption{Example 4. An illustration of the exact coefficient $k$ and the reconstructed coefficients $\tilde{k}$ under fully and partially observed settings, plotted along the diagonal line $y = x$.}
    \label{fig:discontinuous_k_slice}
\end{figure}
Figure~\ref{fig:discontinuous_k_slice} shows a slice of the exact coefficient $k$ and the reconstructed coefficients $\tilde{k}$ under the fully and partially observed settings along the diagonal line $y = x$.
The solid orange line represents the exact coefficient, while the blue and red dashed lines correspond to the reconstructions from the fully and partially observed settings, respectively.
Although slight differences appear near the discontinuities, both reconstructions closely follow the exact coefficient with minimal deviation.

\begin{table}[!ht]
    \centering
    \begin{tabular}{|c|c|c|c|c|c|}
    \hline
    \textbf{Case} & \textbf{Observed} & \textbf{MSE} & \textbf{MAE} & \textbf{PSNR} & \textbf{SSIM} \\ \hline
    \multirow{2}{*}{Discontinuous} 
        & Fully     & $4.36 \times 10^{-4}$ & $2.81 \times 10^{-3}$ & 32.69 & 0.9733 \\ \cline{2-6}
        & Partially & $5.02 \times 10^{-4}$ & $3.30 \times 10^{-3}$ & 32.08 & 0.9676 \\ \hline
    \end{tabular}
    \caption{Example 4. Comparison of MSE, MAE, PSNR, and SSIM for the discontinuous case under fully and partially observed settings.}
    \label{tab:discontinuous_metrics}
\end{table}
Table~\ref{tab:discontinuous_metrics} presents the performance metrics—MSE, MAE, PSNR, and SSIM—for the discontinuous case under both the fully and partially observed settings.
All metrics show slightly better performance in the fully observed setting, indicating a more accurate reconstruction overall.
However, the partially observed setting also achieves highly competitive results.
These findings suggest that the proposed method remains robust even with limited boundary data, particularly since reconstruction errors are primarily localized near discontinuities.

\section*{Acknowledgement}
The work of D. Park was supported by the BK21 FOUR (Fostering Outstanding Universities for Research, No. 202507810000), funded by the Ministry of Education (MOE, Korea) and the National Research Foundation of Korea (NRF-2022R1C1C1003464 and RS-2023-00217116).
The work of S. Lee was supported by the U.S. Department of Energy, Office of Science, Energy Earthshots Initiatives under Award DE-SC-0024703.
The work of S. Moon was supported by the National Research Foundation of Korea (NRF-2022R1C1C1003464 and RS-2023-00217116).

\bibliographystyle{elsarticle-harv}
\bibliography{reference}

\end{document}